\documentclass[final,3p,authoryear]{elsarticle}

\usepackage[T1]{fontenc}
\usepackage[latin9]{inputenc}
\usepackage{array}
\usepackage{multirow}
\usepackage{amsmath}
\usepackage{natbib}
\usepackage{comment}
\usepackage{amsmath}
\usepackage{geometry}
\usepackage{changepage}

\usepackage{adjustbox}

\usepackage{hyperref}

\usepackage{float}
\usepackage{graphicx}

\providecommand{\tabularnewline}{\\}
\floatstyle{ruled}
\newfloat{algorithm}{tbp}{loa}
\providecommand{\algorithmname}{Algorithm}
\floatname{algorithm}{\protect\algorithmname}

\usepackage[british]{babel}

\journal{Not decided}

\begin{document}

\begin{frontmatter}

\title{Truck pooling and scheduling in post-distribution cross-docking with JIT demands and synchronized interchangeability}

\author[label1]{Seyed-Esmaeil Moussavi}
\ead{sesmaeil.moussavi@univ-artois.fr}
\address[label1]{Universit\'{e} d'Artois, Lab-LGI2A, F-62400 B\'{e}thune, France}

\author[label2]{Rahimeh Neamatian Monemi}
\address[label2]{Universit\'{e} de Lille, CRISTAL, F-59655 Lille, France}
\ead{Rahimeh.Monemi@univ-lille.fr}

\author[label3]{Shahin Gelareh}
\ead{shahin.gelareh@univ-artois.fr}
\address[label3]{Universit\'{e} d'Artois, D\'{e}partement R{\&}T, IUT de B\'{e}thune, F-62400 B\'{e}thune, France}




\begin{abstract}

Various operational optimization problems arise in cross-dock synchronization. The combination of the truck scheduling with other decision problems in cross-docking has been targeted by researchers in the recent years. In most of the truck scheduling researches, the load and the destination of the outbound trucks are predetermined. This paper presents a novel cross-docking optimization problem in which the truck scheduling is combined by two assignment problems: the assignment of received loads and the assignment of destinations to outbound trucks (product-truck-destination allocation). Moreover, a Just-In-Time (JIT) strategy is imposed on the destination demands, whereas in most of the previous researches the time windows are imposed to the trucks rather than the demands. An integrated mathematical model is presented for this cross-docking problem. The mathematical model is reinforced by two categories of symmetry breaking constraints. A matheuristic approach, including two mixed-integer mathematical models, is proposed. An adaptive heuristic is presented to solve the first step of the matheuristic algorithm instead of the mathematical model. The solution approaches are analyzed by the extensive computational experiments. The numerical results show the efficiency of the hybrid matheuristic to solve the real size instances of the studied problem.

\end{abstract}


\end{frontmatter}

\section{Introduction} 

Product distribution strategy has a significant impact on the supply chain performance. According to \cite{apte2000effective}, 30\% of an item price is incurred in the distribution process. Hence, lots of companies are currently trying to develop new distribution strategies to efficiently manage their product flow. Cross-docking, as a logistic policy in the distribution part of supply chain, is implemented to improve transport and delivery activities by consolidating the products. Cross-dock is a distribution facility in which the goods arriving by the inbound trucks are unloaded, sorted, consolidated based on their destination and loaded onto the outbound trucks. Inside the cross-dock, the products can be delivered directly to their related outbound trucks or stored in the temporary storage but not more than 24h. On the one hand, warehouses allow companies to make economies in transportation (by making full truckloads), reduce the effect of demand variability and execute additional operations closer to final customers \citep{dolgui2010supply}. On the other hand, cross-docking has the potential to eliminate the storage and retrieval operation functions of a traditional warehouse. As it provides the best of the warehousing and direct delivery strategy, cross-docking has become a popular distribution strategy in practice \citep{Rijal2019}.

The primary purpose of the cross dock is to enable a consolidation of differently sized shipments with the same destination to full truckloads, so that economies in transportation costs can be realized \citep{boysen2010cross}. Indeed, cross docking can decrease distribution network costs by reducing products' delivery time, material handling, and inventory holding costs \citep{kaboudani2018vehicle, theophilus2019truck, zarandi2016constraint}. In the cross-docking operations, an efficient transshipment process is required where inbound and outbound truckloads are synchronized, so that intermediate storage inside the terminal is kept low and on-time deliveries are ensured \citep{boysen2010cross}.

Various decision problems, arising in cross-docking, have been studied in the literature. They could be classified into strategic, tactical and operational problems. Vehicle scheduling, as an operational level optimization problem, is an important concern in the cross-dock planning which is largely considered by researchers. In the classical cross-docking operations, each product is dedicated to a specific outbound truck; But with interchangeability considerations, two products of the same type can be replaced in the predefined load of an outbound truck. Unlike classical interchangeability, in the studied problem, the loads of the outbound trucks are not predefined. Hence, this article presents a new extension of the product interchangeability in which, the load of the outbound trucks and the destination of each truck must be determined by solving the problem. 

In the cross-docks where the product interchangeability is allowed, called post-distribution cross-docking, any arrived pallet must be attributed to an outbound truck (product-truck allocation) or to a destination (product-destination allocation). The studied truck scheduling problem incorporates both mentioned allocations simultaneously. Therefore, it deals with two synchronized interchangeabilities, called product-truck-destination allocation, in a post distribution cross-docking environment.  The integration of these synchronized assignment problems with both inbound and outbound trucks scheduling is the main contribution of this study. Moreover, in the literature, the time windows are mostly defined for the outbound trucks; But in this research, they are assigned to the customers' orders under a JIT principle. 


\cite{boysen2010scheduling} proved that the truck scheduling in a cross-dock facility with one inbound door and one outbound door is an NP hard problem. This paper studies a truck scheduling with multiple inbound and outbound doors combined with two synchronized assignment problems. It can be concluded that the presented problem is also NP hard. Consequently, a matheuristic approach, based on the decomposition of the formulation, is presented to solve the real size problem instances. Afterward, an adaptive heuristic algorithm is provided to improve the performance of the matheuristic.  

The rest of this paper is organized in five sections as follows. After a brief literature review in section \ref{literature}, the studied problem is presented in section \ref{problem description}. The problem is formalized mathematically in section \ref{math model}. The solution approaches (matheuristic and heuristic methods) are explained in section \ref{solution approaches}. The results of the computational experiments are summarized in section \ref{results}. Finally, the conclusion of the paper and the directions for the future studies are presented in section \ref{conclusion}.




\section{Literature}\label{literature}

Cross-docking, as a product distribution strategy that serves as an intermediate node between suppliers and customers in the distribution networks \citep{mousavi2014location, joo2013scheduling}, has received  increasing attention by the researchers over the last two decades. The decision problems in the cross-docking contain a large domain of researches. From strategic to operational, from design to scheduling, from network to local problems, are decided in the cross-dock optimization domain. These decision problems have been classified from the different point of views in the previous researches presented by \cite{Ladier2016}, \cite{buijs2014synchronization}, \cite{van2012cross} and \cite{boysen2010cross}. In this paper, the operational level problems inside the cross-dock environment are targeted.  

Various operational level problems in the cross-dock planning have been studied by the researchers. Transportation planning, dock-door assignment, Vehicle Routing Problem (VRP) and vehicle scheduling are the main optimization problems in this domain. 
For example, a dock-door assignment problem was studied by \cite{Nassief2018} in which the unloading and loading times were considered. Three new MIP formulations were presented and solved by LP relaxation and column generation algorithm. A transportation problem in the cross-dock network was presented by \cite{Musa2010} in which the direct shipment was possible. The objective was to minimize the total transportation cost by determining the best fleet dispatching and consolidation plans. A vehicle routing problem considering split delivery and time window for suppliers/cross-docks/retailers was presented by \cite{wang2017two}, where Less-than-TrackLoad (LTL) was allowed. Two extended approaches based on the simulated annealing and tabu search algorithms were proposed to minimize the total transportation cost.

In recent years, vehicle scheduling problem and its extensions have received more attention from the researchers. In this way, \cite{schwerdfeger2018just} studied an outbound truck scheduling problem in a specific type of cross-dock network, namely one-to-one.  They seek to minimize the number of trucks for the parts deliveries regarding to the JIT intervals demanded by the plant. A heuristic algorithm based on the decomposition of the formulation is presented and compared with binary search and an exact method. \cite{Boysen2010} presented a truck scheduling problem without temporary storage in a cross-dock of food industry distribution. They studied three objectives separately where they developed a dynamic programming and a simulated annealing algorithm to solve this problem. \cite{Chen2009} aimed to minimize the makespan in a cross-dock system. They considered the cross-dock scheduling as a two machines flow-shop problem in which the variable loading/unloading times depended on the products (jobs). They proposed four heuristic algorithms containing two dynamic algorithms to solve the problem.

Moreover, many studies have been carried out on the combinatorial optimization problems in cross-docking. Dock-door assignment combined with VRP and vehicle scheduling are the frequent combinations in this way. For example, \cite{ENDERER201730} presented a combinatorial problem by integrating dock-door assignment and vehicle routing (DAVRP). The assignment of suppliers to the inbound doors and the routing from the outbound doors to the destinations were targeted with aiming to minimize the handling and transport cost. For the combination of  vehicle scheduling with door assignment, \cite{Wisittipanich2017} presented a truck scheduling problem, comprising truck sequencing and door assignment simultaneously. The objective of their problem was to minimize the total operational time (makespan). They proposed a novel variant of Particle Swarm  Optimization (PSO), namely GLNPSO to solve this problem. \cite{Rijal2019} presented another combinatorial problem in which the truck scheduling and dock-door assignment are solved simultaneously. They analyzed the mixed-mode dock-doors in a U-shape cross-docking problem in which a soft time windows for the truck departures and a temporary storage are considered. The objective was to minimize a combination of the costs, containing:  inside traveling, storage and tardiness penalty. An ALNS algorithm was proposed to solve this problem. In another research, \cite{Kusolpuchong2019} studied a truck scheduling and dock-door assignment problem in a temporary storage allowed cross-dock. They proposed a genetic algorithm to optimize three objectives as travel distance of material handling, makespan and maximum storage level. The heterogeneous trucks (three sizes), LTL and the loading/unloading times were taken into account in their problem.

In various vehicle scheduling studies, a time windows was considered to the departure time of the outbound trucks. For example, in the truck scheduling model presented by  \cite{Serrano2017}, a predefined time windows for the arrival trucks is considered as a soft constraint and the objective was to minimize the penalty of the violating these contracted times. A hard time windows for the outbound trucks was taken into account in the truck scheduling door assignment problem presented by \cite{Molavi2018}. Their model proposes the schedules for both inbound and outbound trucks and minimizes the transport cost and the penalty of lateness of delayed shipments by employing a hybrid heuristic algorithm composed of GA and Variable Neighborhood Search (VNS).

The simultaneous consideration of multi-period planning and departure time windows have been studied in a truck scheduling research where \cite{Ghomi2020} proposed an inbound truck scheduling model on a multi-period horizon of time with a temporary storage and variable arrival times for the inbound trucks. A fixed time windows was predefined for the departure of outbound trucks. In another research, \cite{Shahramfard2019} presented a multi-period truck scheduling and door assignment in a cross-dock with regarding to the material handling constraints and assignment of forklifts to the dock-doors. The time window constraints for both inbound and outbound trucks are considered in their schedules calculations. The waiting time for the outbound trucks together with the delayed outbound trucks and holding cost are considered as the objectives and four multi-objective meta-heuristic algorithms are proposed and compared.

The product-truck allocation is another decision problem in the cross-docking domain which arises where the product interchangeability is possible. That means each outbound truck has a list of products and its can be provided by every inbound truck (supplier) under a post-distribution strategy. \cite{tang2010pre} and \cite{yan2009pre} compared the cross-docking operations under post-distribution and pre-distribution considerations. Only a few papers have studied the product-truck allocation (post-distribution considerations) in their truck scheduling. The studies presented by \cite{tootkaleh2016cross}, \cite{assadi2016differential},  \cite{liao2013simultaneous} and \cite{lee2012genetic} considered post-distribution issues in their vehicle planning problems.

In this paper, a post distribution truck scheduling problem is presented in which the load and the destination of the outbound trucks must be determined (product-truck-destination allocation) and  a specific type of time windows (just-in-time) is imposed on the customer demands rather than on the trucks. In order to illustrate the contributions of this study, a comparison with most relevant papers is carried out in table \ref{tab:literature}. All of the mentioned papers incorporated an assignment problem in their post-distribution truck scheduling problems (whether product to truck or product to destination); whereas in this study two interrelated assignments are integrated and synchronized.

\begin{table}

\caption{\label{tab:literature}Related literature on post-distribution truck scheduling in cross-docking}
\begin{adjustbox}{max width=\textwidth}
\begin{tabular}{ccccccccc}
\hline 
Paper & Product allocation & Arrival time & Departure condition & Split d. & LTL & T. storage & Objective & Solution approach\tabularnewline
\hline 
\cite{assadi2016differential} & Product-Truck & Both trucks & JIT-Trucks & No & No & Infinite & Earliness-Tardiness & DE-SA\tabularnewline
\cite{bellanger2013three} & Product-Dest. & Zero & No & No & No & Infinite & Makespan & B\&B\tabularnewline
\cite{ladier2018crossdock} & Product-Dest. & Both trucks & TW-Destinations & Yes & No & Infinite & Inventory Level & TS-IP Decomposition\tabularnewline
\cite{larbi2011scheduling} & Product-Dest. & Inbound & No & Yes & Yes & Infinite & Handling cost & Stochastic heuristic\tabularnewline
\cite{lee2012genetic} & Product-Truck & Zero & No & Yes & Yes & Infinite & Throughput & GA\tabularnewline
\cite{liao2020integrated} & Product-Truck & Zero & TW-Destinations & No & No & Infinite & Earliness-Tardiness & ILS-Greedy search\tabularnewline
\cite{liao2013simultaneous} & Product-Truck & Zero & JIT-Trucks & No & No & Infinite & Tardiness & ACO-Hybrid DE\tabularnewline
\cite{nasiri2018incorporating} & Product-Truck & Zero & TW-Trucks & Yes & No & Infinite & Delivery costs & TSSA\tabularnewline
\cite{Serrano2017} & Product-Dest. & Inbound & No & No & No & Limited & Arrival times & MILP-CPLEX\tabularnewline
\cite{shahmardan2020truck} & Product-Truck & Zero & No & No & Yes & infinite & Makespan & LR-SA\tabularnewline
\cite{shakeri2012robust} & Product-Truck & Zero & Non & No & No & Limited & Makespan & 2phase heuristic\tabularnewline
\cite{tootkaleh2016cross} & Product-Truck & Zero & TW-Trucks & No & No & Infinite & Inventory cost & Constructive heuristic\tabularnewline
\cite{Wisittipanich2017} & Product-Truck & Zero & No & No & No & Infinite & Makespan & Adaptive PSO\tabularnewline
\textbf{Our paper} & \textbf{Product-Truck-Dest.} & \textbf{Inbound} & \textbf{JIT-Demands} & \textbf{Yes} & \textbf{Yes} & \textbf{Limited} & \textbf{Load waiting} & \textbf{Hybrid matheuristic}\tabularnewline
\hline 
\end{tabular}
\end{adjustbox}
\end{table}

\section{Problem description}\label{problem description}

A set of inbound trucks containing the pallets of different products come from the suppliers and arrive to the cross-dock in the estimated times. Under a JIT strategy, each customer needs a number of these products in a specific time windows. A set of outbound trucks must carry the needed products to the customers during the relevant intervals of time. The number of pallets to be transported from an inbound truck (supplier) to an outbound truck (destination) are not predefined; Therefore, the outbound trucks can be fulfilled through any arrived pallets. In other words, the product interchangeability is allowed. The objective is to minimize the total storage time of the pallets inside the cross-dock. The best solution for a pallet is to be unloaded from its related inbound truck and loaded to its allocated outbound truck at the same period. For that purpose, the inbound and outbound trucks must be docked simultaneously. In this case, the waiting time for the pallet will be zero.

The time windows of the customer demands must be respected. Note that, the time windows are assigned to the products rather than to the outbound trucks. Therefore, all of the pallets assigned to an outbound truck have the same due date, which determines the leaving time of that truck.

Furthermore, the outbound trucks can leave the cross-dock with less-than-truckload (LTL) if there is no more pallet available. In this study, the split shipment is possible, where a destination can be served by different trucks. The problem is simplified by some assumptions:

\begin{itemize}

\item The containers which are loaded on a specific outbound truck belong to the same destination.

\item The truck and the pallets are homogeneous.  
\item The un/loading and transshipment times inside the cross-dock are independent of the load size and product type and is the same for all trucks and products.
\item The number of outbound trucks is limited by the constraints, and the transportation time between the cross-dock and every destination is ignored.
\item A working day is composed of a number of time periods.
\item The time during which a truck is docked to the cross-dock (loading/unloading time) in one period. 

\subsubsection*{Interchangeability conditions:}

\item The load of the outbound trucks can be provided by every inbound trucks (suppliers).

\item Each outbound truck can be attributed to any destination

\end{itemize}

The problem is to determine two schedules for the inbound and outbound trucks, and the number of pallets to be transported from each inbound truck to every outbound truck. In other words, in addition to the truck scheduling, the assignment of the suppliers to the destinations for each product type on every outbound trucks is determined. It signifies two synchronized allocations. The allocation of the received products to the outbound trucks and the allocation of the outbound trucks to the destinations.

\begin{figure}

\centering
\includegraphics[width=1\textwidth]{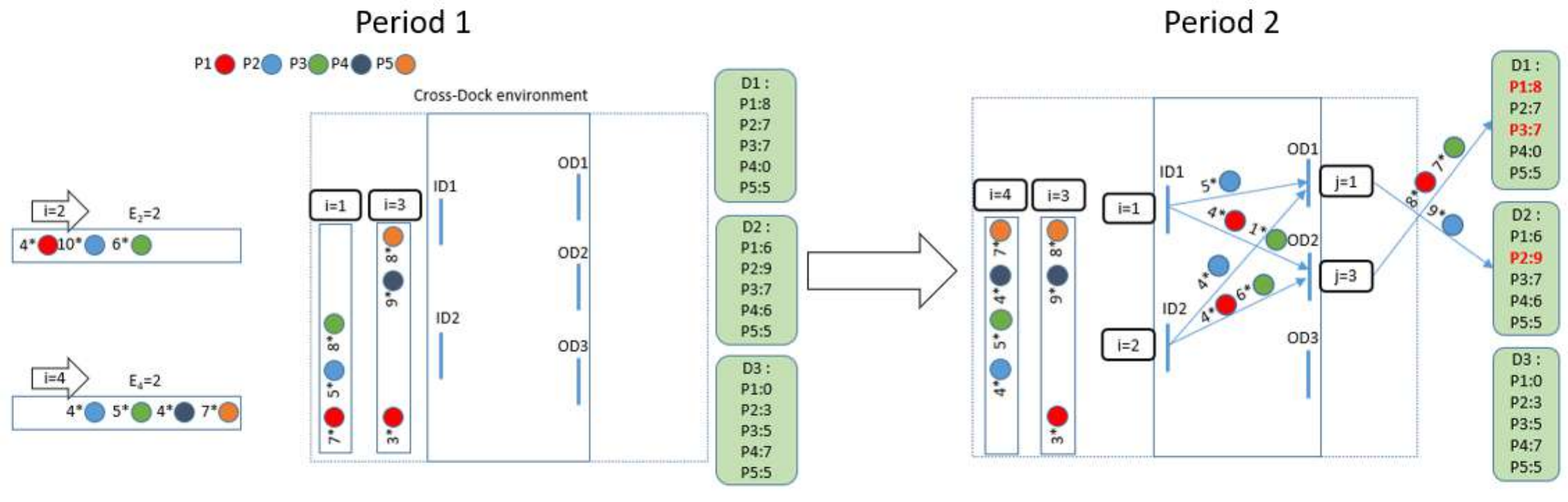}
\caption{Problem description example}
\label{Fig-Problem}
\end{figure}

Figure \ref{Fig-Problem} provides an example for the planning of two periods. Five product types are shown by different colors ($P_{1}$, $P_{2}$,..., $P_{5}$). Two inbound trucks ($i_{1}$ and $i_{3}$) arrive at the first period but they are not docked to the cross-dock. There are three destinations ($D_{1}$, $D_{2}$, $D_{3}$), and each one needs a list of the products under a JIT strategy. In the first period no demand must be delivered; But in the second period, 8 pallets of the first product and 7 pallets of the third product must be delivered to the first destination and 9 pallets of the second product must be shipped to the second destination. In order to meet these demands, two inbound trucks ($i_{1}$ and $i_{2}$) are docked to the inbound doors and unloaded. At the same time, two outbound trucks ($j_{1}$ and $j_{3}$) are docked to the outbound side and pooled by the received pallets. The first outbound truck contains 9 pallets of product 2 and that is assigned to destination 2. The transshipment between every pair of trucks are shown in the figure for this small example. 

The best solution for a pallet is to be unloaded from its related inbound truck, and loaded to its allocated outbound truck at the same period. For this purpose, both inbound and outbound trucks must be docked simultaneously. In this case, the waiting time of the pallet will be zero. 

An integrated mathematical model, which involves truck scheduling and product-truck-destination allocation, is presented for this post-distribution cross-docking problem.

\subsection{Mathematical model}\label{math model}

A mixed-integer mathematical model is proposed for the aforementioned problem. The objective is to minimize the waiting time of the pallets inside the cross-dock. It means the time interval length between unloading a pallet from inbound truck and loading it on the outbound truck. The parameters, variables, and problem formulation are presented as follows.

\subsubsection{Parameters and variables}

The sets and the parameters which are used to formulate the problem mathematically are presented in table \ref{table param}. The integer and binary variables, which are employed to model the problem, are listed in table \ref{table var}.

\begin{flushright}
\begin{table}[h]
\caption{Parameter and set notation}
\label{table param}
\small
\begin{adjustbox}{max width=\textwidth}
\begin{tabular}{ll}
\hline 
\textbf{Sets} & \tabularnewline
$I$ & Set of inbound trucks; $I=\{1,2,...,m\}$;\tabularnewline
$J$ & Set of outbound trucks; $J=\{1,2,...,n\}$;\tabularnewline
$P$ & Set of product types; $P=\{1,2,...,k\}$;\tabularnewline
$D$ & Set of destinations; $D=\{1,2,...,f\}$;\tabularnewline
$T$ & Set of time periods; $T=\{1,2,...,r\}$;\tabularnewline
\textbf{Parameters} & \tabularnewline
$m$ & Total number of inbound trucks; \tabularnewline
$n$ & Total number of outbound trucks; \tabularnewline
$k$ & Total number of product types; \tabularnewline
$f$ & Total number of destinations (customers); \tabularnewline
$r$ & Total number of time periods during a working day; \tabularnewline
$ID$ & Total number of inbound doors; \tabularnewline
$OD$ & Total number of outbound doors;\tabularnewline
$C$ & Capacity of the trucks; (homogeneous trucks);\tabularnewline
$E_{i}$ & Arriving time of the inbound trucks;\tabularnewline
$L_{ip}$ & Total number of pallets of product $p$ in inbound truck $i$;\tabularnewline
$R_{pdt}$ & Number of pallets $p$ needed by destination $d$ in period $t$ (JIT customer demands);\tabularnewline
$PC$ & Penalty cost for each non-delivered pallet of a customer order;\tabularnewline
\multirow{1}{*}{$RB_{pdt}$} & $\begin{cases}
=1,\text{ if customer d needs product p in period t;}\\
=0,\text{ otherwise.}
\end{cases}$\tabularnewline
\hline 
\end{tabular}
\end{adjustbox}
\end{table}
\par\end{flushright}

\begin{flushleft}
\begin{table}[h]


\caption{Decision variables}
\label{table var}
\small
\begin{adjustbox}{max width=\textwidth}
\begin{tabular}{ll}
\hline 
$S_{ijpd}$ & Number of pallets $p$ delivered to destination $d$ by outbound truck $j$ which is supplied by inbound truck $i$\tabularnewline
\multirow{1}{*}{$SB_{ijt}$} & $\begin{cases}
=1,\text{ if there is a shipment from inbound truck i to outbound truck j at time period t; }\\
=0,\text{ otherwise;}
\end{cases}$\tabularnewline
\multirow{1}{*}{$y_{jt}$} & $\begin{cases}
=1\text{ if outbound truck \ensuremath{j} is docked to the cross-dock (is loaded) in period \ensuremath{t}; }\\
=0,\text{ otherwise;}
\end{cases}$\tabularnewline
\multirow{1}{*}{$h_{it}$} & $\begin{cases}
=1\text{ if inbound truck \ensuremath{i} is docked to the cross-dock (is unloaded) in period \ensuremath{t}; }\\
=0,\text{ otherwise;}
\end{cases}$\tabularnewline
\multirow{1}{*}{$q_{jd}$} & $\begin{cases}
=1\text{ if outbound truck \ensuremath{j} is allocated to destination \ensuremath{d}; }\\
=0,\text{ otherwise;}
\end{cases}$\tabularnewline
\multirow{1}{*}{$WB_{jpd}$} & $\begin{cases}
=1,\text{ if outbound truck \ensuremath{j} carries pallets \ensuremath{p} to destination \ensuremath{d};}\\
=0,\text{ otherwise;}
\end{cases}$\tabularnewline
$V_{pd}$ & Number of pallets $p$ which are requested and not supplied to destination $d$; (uncovered demands)\tabularnewline
$St_{t}$ & Number of pallets remained inside the cross-dock at the end of period $t$;\tabularnewline
\hline 
\end{tabular}
\end{adjustbox}
\end{table}

\end{flushleft}

\subsubsection{Problem formulation}

The objective function and the constraints of the studied problem is formulated as follows.

\begin{equation}
Minimize\:Z=\sum_{i\in I}\sum_{j\in J}(((\sum_{t\in T}t\times y_{jt})\times (\sum_{t\in T}SB_{ijt}))-((\sum_{t\in T}t\times h_{it})\times (\sum_{t\in T}SB_{ijt})))+PC\times\sum_{p\in P}\sum_{d\in D}V_{pd}\label{eq:1}
\end{equation}

$Subject\:to:$

\begin{align}	
&\sum_{d\in D}\sum_{j\in J}S_{ijpd}\leq L_{ip} \qquad &\forall i\in I,\:p\in P\label{eq:2}\\
&\sum_{i\in I}\sum_{p\in P}\sum_{d\in D}S_{ijpd}\leq C \qquad &\forall j\in J\label{eq:4}\\
&\sum_{i\in I}\sum_{j\in J}S_{ijpd}+V_{pd}=\sum_{t\in T}R_{pdt}\qquad &\forall d\in D,\:p\in P\label{eq:5}\\
&(\sum_{p\in P}\sum_{d\in D}S_{ijpd})/\sum_{p\in P}L_{ip}\leq\sum_{t\in T}SB_{ijt}\leq\sum_{p\in P}\sum_{d\in D}S_{ijpd} \qquad &\forall i\in I,\:j\in J\label{eq:6}\\
&\sum_{t\in T}SB_{ijt}\leq1\qquad &\forall i\in I,\:j\in J\label{eq:8}\\
&\sum_{i\in I}h_{it}\leq ID\qquad &\forall t\in T\label{eq:9}\\
&\sum_{t\in T}h_{it}=1\qquad &\forall i\in I\label{eq:10}\\
&\sum_{j\in J}y_{jt}\leq OD\qquad &\forall t\in T\label{eq:12}\\
&\sum_{t\in T}y_{jt}\leq1\qquad &\forall j\in J\label{eq:13}\\
&\sum_{t\in T}t\times h_{it}\leq (\sum_{t\in T}t\times SB_{ijt})+(r\times(1-\sum_{t\in T}SB_{ijt}))\qquad &\forall i\in I,\:j\in J\label{eq:17}\\
&\sum_{t\in T}t\times SB_{ijt}\leq \sum_{t\in T}t\times y_{jt}\qquad &\forall i\in I,\:j\in J\label{eq:18}\\
&(\sum_{i\in I}S_{ijpd})/C \leq WB_{jpd}\leq \sum_{i\in I}S_{ijpd}\qquad &\forall j\in J,\:p\in P,\:d\in D\label{eq:20}\\
&\sum_{p\in P}WB_{jpd}/r\leq q_{jd}\qquad &\forall j\in J,\:d\in D\label{eq:30}\\
&\sum_{d\in D}q_{jd}\leq1\qquad &\forall j\in J\label{eq:31}\\
&(\sum_{t\in T}t\times y_{jt})\times WB_{jpd}=(\sum_{t\in T}t\times RB_{pdt})\times WB_{jpd}\qquad &\forall j\in J,\:p\in P,\:d\in D\label{eq:22}\\
&Capacity \: constraints:& \nonumber \\
&St_{t}=St_{t-1}+(\sum_{i\in I}(h_{it}\times\sum_{p\in P}L_{ip})-\sum_{j\in J}(y_{jt}\times\sum_{p\in P}\sum_{d\in D}\sum_{i\in I}S_{ijpd}))\qquad &\forall t\in T \; | \; St_{0}=0\
\label{eq:36}\\
&St_{t}\leq ID\times C\qquad &\forall t\in T\label{eq:37}
\end{align}

The objective function (equation \eqref{eq:1}) minimizes the waiting time of the pallets which must be shipped from inbound trucks to the outbound trucks inside the cross-dock, while minimizing the undelivered demands by adding a penalty cost per pallet. The formulation is nonlinear which must be linearized. For this purpose, the equation is divided to three parts among which two first parts are nonlinear.

The first part ($(\sum_{t\in T}t\times y_{jt})\times (\sum_{t\in T}SB_{ijt})$) can be defined as a set of integer variables ($qy_{ij}$), indicating the docking period of the outbound truck $j$, if there is any shipment from inbound truck $i$ to outbound truck $j$. If there is no shipment between two trucks, $qy_{ij}$ equals zero. Similarly for the inbound trucks, the second part of the objective function ($(\sum_{t\in T}t\times h_{it})\times (\sum_{t\in T}SB_{ijt})$) indicates the docking period of inbound truck $i$, if this truck contains the products to be shipped to outbound truck $j$. This part of the equation can be defined as another set of integer variables ($qh_{ij}$). 

Aforementioned nonlinear parts of the objective function ($qy_{ij}=\sum_{t\in T}t\times y_{jt})\times (\sum_{t\in T}SB_{ijt})$ and $qh_{ij}=\sum_{t\in T}t\times h_{it})\times (\sum_{t\in T}SB_{ijt})$) are linearized by using two sets of inequalities \eqref{eq:28} and \eqref{eq:26}, respectively.

\begin{equation}
\begin{cases}
qy_{ij}\leq r\times (\sum_{t\in T}SB_{ijt})\qquad & \forall i\in I,\:j\in J\\
(\sum_{t\in T}t\times y_{jt})-qy_{ij}\leq r\times(1-(\sum_{t\in T}SB_{ijt}))\qquad & \forall i\in I,\:j\in J\\
(\sum_{t\in T}t\times y_{jt})-qy_{ij}\geq E_{i}-(1-(\sum_{t\in T}SB_{ijt}))\qquad & \forall i\in I,\:j\in J
\end{cases}\label{eq:28}
\end{equation}

\begin{equation}
\begin{cases}
qh_{ij}\leq r\times (\sum_{t\in T}SB_{ijt})\qquad & \forall i\in I,\:j\in J\\
(\sum_{t\in T}t\times h_{it})-qh_{ij}\leq r\times(1-(\sum_{t\in T}SB_{ijt}))\qquad & \forall i\in I,\:j\in J\\
(\sum_{t\in T}t\times h_{it})-qh_{ij}\geq E_{i}-(1-(\sum_{t\in T}SB_{ijt}))\qquad & \forall i\in I,\:j\in J
\end{cases}\label{eq:26}
\end{equation}

Constraints \eqref{eq:2} guarantee that the sum of the products delivered to the customers does not exceed the products carried to the cross-dock. Constraints \eqref{eq:4} limit the load of each outbound truck by its capacity. Constraints \eqref{eq:5} express
that the load delivered to a customer equals the load that customer requires, except in cases in which supplying the whole order is not feasible. Inequalities \eqref{eq:6} define a set of binary decision variables ($SB_{ijt}$) indicating that if $\sum_{p\in P}\sum_{d\in D}S_{ijpd}>0$, it means there is a shipment from inbound truck $i$ to outbound truck $j$, hence $\sum_{t\in T}SB_{ijt}=1$, whereas if $\sum_{p\in P}\sum_{d\in D}S_{ijpd}=0$, hence $\sum_{t\in T}SB_{ijt}=0$.  Constraints \eqref{eq:8} ensure that the shipment from an inbound truck to an outbound truck (if there is any) is occurred in one period of time.

Constraints \eqref{eq:9} guarantee that the number of inbound trucks which are docked simultaneously does not exceed the number of inbound doors. Equation \eqref{eq:10} signifies that all of the inbound trucks must be docked during a working-day. Similarly, constraints \eqref{eq:12} limit the number of outbound trucks which are docked at each period by the number of outbound doors. Constraints \eqref{eq:13} show that it is not necessary that all of the outbound trucks are docked and they are used if needed. Constraints \eqref{eq:17} guarantee that, if there is any shipment from inbound truck $i$ to outbound truck $j$, this shipment occurs after docking time of the inbound truck. Constraints \eqref{eq:18} guarantee that all of the shipments to an outbound truck is occurred before or during the docking period of the outbound truck. 

Constraints \eqref{eq:20} define another set of binary variables showing whether outbound truck $j$ transports product $p$ to destination $d$, or not. Constraints \eqref{eq:30} ensure that an outbound truck can transport the pallets only to its attributed destination. Constraints \eqref{eq:31} limit the outbound trucks to visit only one destination. 

Under a JIT strategy, constraints \eqref{eq:22} indicate that an outbound truck is pooled (docked) in the period, which is imposed by the customer order. The right side of this equation is non-linear which can be defined as  a set of variables ($WY_{jpd}$) announcing the docking period of outbound truck $j$, if this truck transports product $p$ to destination $d$, otherwise, it equals zero ($WY_{jpd}=(\sum_{t\in T}t\times y_{jt})\times WB_{jpd}$). This non-linear equation is linearized by employing the following inequalities:

\begin{equation}
\begin{cases}
WY_{jpd}\leq r\times WB_{jpd}\qquad & \forall j\in J,\:p\in P,\:d\in D\\
(\sum_{t\in T}t\times y_{jt})-WY_{jpd}\leq r\times(1-WB_{jpd})\qquad & \forall j\in J,\:p\in P,\:d\in D\\
(\sum_{t\in T}t\times y_{jt})-WY_{jpd}\geq0\qquad & \forall j\in J,\:p\in P,\:d\in D
\end{cases}\label{eq:23}
\end{equation}

\subsubsection*{Capacity constraints:}

At the end of each period, the total number of pallets remained in the temporary storage is obtained by equation \eqref{eq:36}. The temporary storage in period $t$ is calculated by adding the total number of pallets unloaded from the inbound trucks, to the temporary storage of the precedent period and subtracting the total number of pallets loaded onto the outbound trucks. In equation \eqref{eq:36}, the third part of the right side of the formulation ($\sum_{j\in J}(y_{jt}\times\sum_{p\in P}\sum_{d\in D}\sum_{i\in I}S_{ijpd})$) is nonlinear. For the linearization, a set of integer variables ($LJ_{jt}$) are defined which are equivalent to this part of the equation ($LJ_{jt}=y_{jt}\times\sum_{p\in P}\sum_{d\in D}\sum_{i\in I}S_{ijpd}$). $LJ_{jt}$ and present the total number of pallets loaded onto outbound truck $j$ in period $t$. It's related nonlinear formulation is linearized by adding the following inequalities:

\begin{equation}
\begin{cases}
LJ_{jt}\leq C\times y_{jt}\qquad & \forall j\in J,\:t\in T\\
\sum_{p\in P}\sum_{d\in D}\sum_{i\in I}S_{ijpd}-LJ_{jt}\leq C\times(1-y_{jt})\qquad & \forall j\in J,\:t\in T\\
\sum_{p\in P}\sum_{d\in D}\sum_{i\in I}S_{ijpd}-LJ_{jt}\geq0\qquad & \forall j\in J,\:t\in T
\end{cases}\label{eq:33}
\end{equation}

Constraints \eqref{eq:37} limit the temporary storage level that can not exceed a constant value which is calculated by multiplying the number of inbound doors by the capacity of a truck.

After formalizing the problem mathematically, a solution procedure is proposed to this problem. The steps and the algorithms of the solution procedure are explained in the next section.

\section{Solution method}\label{solution approaches}

In order to solve the studied problem, firstly, a set of symmetry breaking constraints are added to the mathematical model to make that more efficient. At the second step, a two-step matheuristic method is presented which is based on the decomposition of the formulation. As the first step of the matheuristic is still complex, a adaptive heuristic is proposed as an alternative to solve this step of matheuristic. Hence, a hybrid matheuristic, composed of an adaptive heuristic (for the first step) and a mathematical model for the second step, is presented to solve the studied combinatorial optimization problem.

\subsection{Symmetry breaking}

In any solution of the above-mentioned mathematical model, three attributes are determined for each outbound truck as: the load, the docking period and the destination. Solution 1 is a part of an entire solution of the problem in which $Y_{j}$ and $D_{j}$ imply respectively the docking period and the assigned destination of outbound truck $j$.

$Solution\:1\::\:\begin{cases}
Y_{1}=A_{1},...,Y_{i}=A_{i},...,Y_{j}=A_{j},...,Y_{n}=A_{n};\\
D_{1}=B_{1},...,D_{i}=B_{i},...,D_{j}=B_{j},...,D_{n}=B_{n};
\end{cases}$

Since the destinations and the loads of the outbound trucks are not predefined, thus another solution could be easily produced by simultaneously substituting the docking periods and the destinations of two trucks. As an example, solution 2 shows a part of a solution derived from solution 1 and obtained by replacing the docking period and the destination of truck $i$ by the ones of truck $j$.

$Solution\:2\::\:\begin{cases}
Y_{1}=A_{1},...,Y_{i}=A_{j},...,Y_{j}=A_{i},...,Y_{n}=A_{n};\\
D_{1}=B_{1},...,D_{i}=B_{j},...,D_{j}=B_{i},...,D_{n}=B_{n};
\end{cases}$

Because the outbound trucks are homogeneous, i.e. they have the same capacity, cost, entering time, availability and speed, ..., solution 2 (and any other solution obtained in the same way) results in the same objective value as solution 1. Hence, in each solution node of the problem, a group of symmetric solutions is raised that causes an exceedingly vast branch and bound three. Consequently, by keeping the values of other variables, every permutation of the outbound trucks ($n!$) over theirs docking periods and destinations simultaneously, gives a symmetric solution. 

In order to avoid or decrease the symmetric solutions, appropriate inequalities called ``symmetry breaking constraints'' could be imposed to the model. In this paper, two categories of the symmetry breaking inequalities are proposed to take away all or a large number of symmetric solutions. Considering an ordering strategy to assign the trucks to the cross-dock is the first issue which is taken into account to handle the symmetries. In this way, inequality \eqref{eq:SB1} is imposed to the problem, knowing that it does not touch the solutions nodes and the feasible space. 

\begin{equation}
Y_{1}\leq Y_{2}\leq...\leq Y_{n} \label{eq:SB1}
\end{equation}

This inequality keeps only one or a small number of solutions at each node. It signifies that outbound truck $i$ is docked to the cross-dock before or at the same period as the outbound truck $j$, if $i<j\;(i,j\in J)$.

By employing these constraints, a large number of symmetric solutions are removed, but there still remains some sub-symmetries inside each planning period. For $n$ outbound trucks, $n!$ possible permutations signifies that in each node, there exists $n!$ ($\prod_{1\leq i\leq n}i$) different solutions; But by applying inequality \eqref{eq:SB1}, this number is reduced to : $n'_{1}!\times n'_{2}!\times...\times n'_{r}!$, in which, $n'_{t}$ implies the number of outbound trucks assigned to time period $t$ $(t\in T)$, such that $n'_{1}+n'_{2}+...+n'_{r}=n$.

\textbf{\textit{Lemma 1}} : $If\:(n'_{1}+n'_{2}+...+n'_{r}=n)$ $then:$ $(\prod_{1\leq i\leq n'_{1}}i)\times(\prod_{1\leq i\leq n'_{2}}i)\times...\times(\prod_{1\leq i\leq n'_{r}}i)\leq(\prod_{1\leq i\leq n}i)$.

\textbf{\textit{Lemma 2}} : $If\:(n>0)$ $then:(\prod_{1\leq i\leq n'_{1}}i)\times(\prod_{1\leq i\leq n'_{2}}i)\times...\times(\prod_{1\leq i\leq n'_{r}}i)>0$.

Two above-mentioned lemmas signify that the number of symmetric solutions is decreased while keeping at least one solution in each node, i.e. these inequalities do not touch the nodes and feasible space of the problem.

The second category of Symmetry Breaking Constraints (SBC) is to avoid the sub-symmetries inside each planning period and it concerns the assignment of the trucks to the destinations. Even by taking into account the inequality \eqref{eq:SB1}, each solution node still contains multiple sub-symmetry groups. Inside a planning period, any permutation of the assigned trucks $(J'_{t}\subset J)$ over their destinations results in the same solution node and objective value. Similar to the first category of SBC, an ordering system is defined to keep one or a small number of solutions of each node and remove the others. This ordering inequality is defined as follows: 

\begin{equation}
D_{i}<D_{j}\qquad\forall t\in T\;and\;i,j\in\{J'_{t}\mid i<j\};\label{eq:SB2}
\end{equation}

This inequality signifies that among the trucks which are docked simultaneously, the truck with smaller number is assigned to the destination with smaller number. 

Under mixed-integer linear programming conditions, two mentioned SB inequalities are formulated as follows: 

\begin{equation}
Y_{j}-Y_{j-1}\geq0\qquad\forall j\in\{J|\;j\geq1\}\label{eq:38}
\end{equation}

\begin{equation}
(\sum_{d\in D}d\times q_{jd})\times y_{jt}\leq((\sum_{d\in D}d\times q_{gd})\times y_{gt})+f\times(1-y_{gt})\qquad\forall j,g\in \{J|\;j<g\},\:t\in T\label{eq:42}
\end{equation}

Constraints \eqref{eq:38} guarantee that the outbound trucks are assigned to the cross-dock by increasing order of their numbers. Constraints \eqref{eq:42} define the assignment of the trucks to the destinations for the trucks which are docked at the same period. By these constraints, the smaller truck number is assigned to the smaller destination number. The left side of this inequality is nonlinear which must be linearized. For this purpose, a set of integer variables ($DT_{jt}$) is defined which is equivalent to the left side of inequality \eqref{eq:42} :  $DT_{jt}=(\sum_{d\in D}d\times q_{jd})\times y_{jt}$. This equation is linearized by using inequalities \eqref{eq:41} as the follows:

\begin{equation}
\begin{cases}
DT_{jt}\leq f\times y_{jt}\qquad & \forall j\in J,\:t\in T\\
\sum_{d\in D}d\times q_{jd}-DT_{jt}\leq f\times(1-y_{jt})\qquad & \forall j\in J,\:t\in T\\
\sum_{d\in D}d\times q_{jd}-DT_{jt}\geq0\qquad & \forall j\in J,\:t\in T
\end{cases}\label{eq:41}
\end{equation}

In the same manner, the first part of the right side of inequality \eqref{eq:42} ($(\sum_{d\in D}d\times q_{gd})\times y_{gt}$) is linearized by considering the variables $DT_{gt}$.

By adding SBCs (equations \eqref{eq:38} and \eqref{eq:42}), the efficiency of the model is significantly improved; But in order to solve the real size problem instances, the more efficient solution approaches are needed. Regarding to the complexity of the model, the matheuristic method seems to be a suitable method to solve the studied problem. Hence, a two-step matheuristic algorithm is proposed and explained as follows.

\subsection{Matheuristic method}

A matheuristic method based on the decomposition of the formulation has been applied to solve the studied combinatorial optimization problem. By this method, the complex problem is decomposed to two sub-problems where the first sub-problem is solved and it's decision variables are used as the parameters to the second sub-problem. The complexity of the problem has been significantly decreased by using this decomposition. For the studied problem, the scheduling of the outbound trucks, the assignment of the outbound trucks to the destinations and the load (number of pallets of each product) of the outbound trucks is solved at the first step. Actually, all the variables related to the outbound of cross-dock and the destinations are considered in the first sub-problem. 

In the second step, load and dock time of the outbound trucks are used as the parameters to schedule the inbound trucks and to determine the number of pallets that must be shipped from the inbound trucks to each outbound truck. Actually, the assignment of the product/supplier to the trucks, and accordingly, to the destinations are determined in this step. All constraints related to the inbound part of cross-dock, the load of the inbound trucks and the shipment times and amounts from the inbound trucks to the outbound trucks are considered in this sub-problem.

\subsubsection{Step 1 : Outbound scheduling and destination assignment}

Most of the variables, which are used in the first sub-problem, are the same as the variables of the main problem defined in the previous section. Moreover, two new sets of variables are applied  to build the connection between the first and the second models. These  variables are defined as follows :

\begin{itemize}
\item $x_{jp}$ : Number of pallets of product $p$ which are placed in outbound truck $j$.
\item $U_{jpdt}$ : Number of pallets of product $p$ which are placed in outbound truck $j$, if this truck is assigned to destination $d$ and is docked (loaded) at the cross-dock in period $t$.
\end{itemize}

Two mentioned sets of variables transform the necessary information from the first model to solve the second model in an optimal manner.  

In this sub-problem, the constraints related to the demands, outbound trucks and doors, product and truck assignment (assignment of the products to the outbound trucks and assignment of these trucks to the destinations) are considered. The objective of this step is to minimize the outstanding demands.

The mathematical model of the first sub-problem is presented as follows :

\begin{equation}
Minimize\:Z1=\sum_{p\in P}\sum_{d\in D}\sum_{t\in T}V_{pdt},\label{eq:301}
\end{equation}

$Subject\:to:$

Assignment of the products and outbound trucks :
\begin{equation}
W_{jpd}= x_{jp}\times q_{jd}\qquad\forall j\in J,\:p\in P,\:d\in D\label{eq:302}
\end{equation}

\begin{equation}
\sum_{d\in D}q_{jd}\leq1\qquad\forall j\in J\label{eq:131}
\end{equation}

Schedule of the outbound trucks :
\begin{equation}
U_{jpdt}=W{jpd}\times y_{jt}\qquad\forall j\in J,\:p\in P,\:d\in D,\:t\in T\label{eq:303}
\end{equation}

\begin{equation}
\sum_{t\in T}y_{jt}\leq1\qquad\forall j\in J\label{eq:113}
\end{equation}

\begin{equation}
\sum_{t\in T}y_{jt}=\sum_{d\in D}q_{jd}\qquad\forall j\in J\label{eq:305}
\end{equation}

Demand/JIT constraints :
\begin{equation}
\sum_{j\in J}U_{jpdt}+V_{pdt}=R_{pdt}\qquad\forall p\in P,\:d\in D,\:t\in T\label{eq:105}
\end{equation}

Outbound doors/Trucks constraints :
\begin{equation}
\sum_{j\in J}y_{jt}\leq OD\qquad\forall t\in T\label{eq:112}
\end{equation}

\begin{equation}
\sum_{p\in P}x_{jp}\leq C\qquad\forall j\in J\label{eq:104}
\end{equation}

\begin{equation}
\sum_{t\in T}\sum_{d\in D}U_{jpdt}=x_{jp}\qquad\forall j\in J,\:p\in P\label{eq:304}
\end{equation}

The objective function of the model (equation \eqref{eq:301}) minimizes the demands which are not satisfied punctually under JIT strategy. That means the total number of pallets of different products which are not delivered at the demanded time to the costumers. Constraints \eqref{eq:302} and \eqref{eq:131} are related to the assignment of the pallets to the outbound trucks and also assignment of the outbound trucks to the destinations. Constraints \eqref{eq:302} define the destination and the load of each trucks and constraints \eqref{eq:131} guarantee that each truck is assigned to maximum one destination.

Constraints \eqref{eq:303}, \eqref{eq:113} and \eqref{eq:305} are concerned with the schedule of outbound trucks. Constraints \eqref{eq:303} define the time at which each truck leaves the cross-dock. Constraints \eqref{eq:113} restrain the outbound truck to be dock to the cross-dock in one and only one period. Constraints \eqref{eq:305} signify that if a truck is assigned to a destination, it must be docked to the cross-dock (must be loaded).

The next category of the constraints are related to the demands under a JIT strategy. Constraints \eqref{eq:105} imply that the outbound trucks ship the pallets to the destinations at the times and in the amounts which are precised and predefined by the customers (destinations); Except for the demands or a part of demands which are not satisfied ($V_{pdt}$). 

Three next sets of constraints (\eqref{eq:112}, \eqref{eq:104} and \eqref{eq:304}) are to define the outbound doors and outbound trucks restrictions. Constraints \eqref{eq:112} signify that the number of trucks docked to the cross-dock at the same times must not exceed the number of outbound doors. Constraints \eqref{eq:104} limit the load (number of pallets) of each outbound truck by it's capacity. Constraints \eqref{eq:304} is to calculate the load (number and type of pallets) which is placed in each outbound truck.

By solving  this step, the value of certain variables are obtained. $x_{jp}$, $q_{jd}$ and $y_{jt}$ are the results of the first model, that show the load, destination and docking period of the outbound trucks.  
 
\subsubsection{Step 2 : inbound scheduling and load assignment}

The value of the variables which are obtained by solving the first model ($x_{jp}$, $q_{jd}$ and $y_{jt}$), are used as the parameters in this step (second sub-problem). The constraints related to the inbound trucks (supply constraints), inbound doors and shipment from inbound trucks to outbound trucks are considered in the second model.

Apart from the main variables, three other sets of variables are defined in this section in order to manipulate (confirm, adjust or modify) the solutions resulted from the first step and obtain the best value for the objective function. These variables are defined as follows :

\begin{itemize}
\item $A_{ijp}$ : (Decision variables) The final number of pallets of product $p$ shipped from inbound truck $i$ to outbound truck $j$ 
\item $B_{ip}$ : Number of pallets of product $p$ from inbound truck $i$ which are not shipped to any outbound truck.
\item $G_{jp}$ : Number of pallets of product $p$ which are assigned to outbound truck $j$ at the first step, but by considering the inbound constraints, they can not be arrived to the outbound truck. 
\end{itemize}

Objective of this sub-problem is the same as the objective of the main problem. This step of the matheuristic is formalized mathematically as follows :

\begin{equation}
Minimize\:Z=\sum_{i\in I}\sum_{j\in J}(((\sum_{t\in T}t\times y_{jt})\times \sum_{t\in T}SB_{ijt})-((\sum_{t\in T}t\times h_{it})\times \sum_{t\in T}SB_{ijt}))+PC\times\sum_{i\in I}\sum_{p\in P}B_{ip}\label{eq:201}
\end{equation}

$Subject\:to:$

Supply constraints :
\begin{equation}
\sum_{j\in J}A_{ijp}+B_{ip}= L_{ip}\qquad\forall i\in I,\:p\in P\label{eq:202}
\end{equation}


\begin{equation}
\sum_{i\in I}A_{ijp}+G_{jp}=x_{jp}\qquad\forall j\in J,\:p\in P\label{eq:205}
\end{equation}


Inbound doors and trucks constraints
\begin{equation}
\sum_{i\in I}h_{it}\leq ID\qquad\forall t\in T\label{eq:209}
\end{equation}

\begin{equation}
\sum_{t\in T}h_{it}=1\qquad\forall i\in I\label{eq:210}
\end{equation}

Schedule of inbound trucks and shipments inside cross-dock :
\begin{equation}
(\sum_{p\in P}A_{ijp})/\sum_{p\in P}L_{ip}\leq\sum_{t\in T}SB_{ijt}\leq\sum_{p\in P}A_{ijp}\qquad\forall i\in I,\:j\in J\label{eq:206}
\end{equation}

\begin{equation}
\sum_{t\in T}SB_{ijt}\leq1\qquad\forall i\in I,\:j\in J\label{eq:208}
\end{equation}

\begin{equation}
\sum_{t\in T}t\times h_{it}\leq (\sum_{t\in T}t\times SB_{ijt})+(r\times(1-\sum_{t\in T}SB_{ijt}))\qquad\forall i\in I,\:j\in J\label{eq:217}
\end{equation}

\begin{equation}
\sum_{t\in T}t\times SB_{ijt}\leq \sum_{t\in T}t\times y_{jt}\qquad\forall i\in I,\:j\in J\label{eq:218}
\end{equation}

Capacity constraints :
\begin{equation}
St_{t}=St_{t-1}+(\sum_{j\in J}(y_{jt}\times\sum_{p\in P}\sum_{i\in I}A_{ijp})-\sum_{i\in I}(h_{it}\times\sum_{p\in P}L_{ip}))\qquad \forall t\in T \; | \; St_{0}=0\
\label{eq:236}
\end{equation}

\begin{equation}
St_{t}\leq ID\times C\qquad\forall t\in T\label{eq:237}
\end{equation}

Similar to the objective of the main model, the objective function of the model aims to minimize the overall non-satisfied  demands by considering both inbound (supply) and outbound constraints while minimizing the waiting times for the unloaded pallets. It means, for each shipment between an inbound and an outbound truck, the time between unloading the pallets from inbound and loading the pallets into the outbound trucks is targeted to be minimized. A part of it's formulation is nonlinear $((\sum_{t\in T}t\times h_{it})\times \sum_{t\in T}SB_{ijt})$, which is linearized as follows :

\begin{equation}
\begin{cases}
qh_{ij}\leq r\times \sum_{t\in T}SB_{ijt}\qquad & \forall i\in I,\:j\in J\\
(\sum_{t\in T}t\times h_{it})-qh_{ij}\leq r\times(1-\sum_{t\in T}SB_{ijt})\qquad & \forall i\in I,\:j\in J\\
(\sum_{t\in T}t\times h_{it})-qh_{ij}\geq E_{i}-(1-\sum_{t\in T}SB_{ijt})\qquad & \forall i\in I,\:j\in J
\end{cases}
\end{equation}

First category of the constraints (\eqref{eq:202} and \eqref{eq:205}) implies the supply restrictions of the cross-dock. Constraints \eqref{eq:202} guarantee that the shipments from an inbound trucks does not exceed the load of the truck. By these constraints, the number of pallets of each product which are not delivered to the outbound trucks (destinations) is calculated. Similarly, constraints \eqref{eq:205} signify that the number of pallets of each product which are shipped to each outbound truck does not exceed the number which is calculated in the first step ($x_{jp}$). The values of the variables $G_{jp}$ is calculated by these constraints. The positive value of $G$ means that it is not possible to response to all of the demands  which are programmed in step one.

Similar to the main model, constraints \eqref{eq:209} and \eqref{eq:210} are related to the number of inbound doors and the time period at which the inbound trucks are docked. Next category of the constraints (\eqref{eq:206} to \eqref{eq:218}) are concerned with the scheduling of the inbound trucks and shipment times inside the cross-dock (between inbound and outbound trucks). Constraints \eqref{eq:206} imply that if there is a shipment from an inbound to an outbound truck, a time period must be assigned to this shipment. Constraints \eqref{eq:208} signify that the shipment from an inbound to an outbound truck, if there is any, is occurred in one and only one period of time. Constraints \eqref{eq:217} guarantee that the shipment from an inbound truck is occurred at the same period or after docking period of the truck. Similarly, constraints \eqref{eq:218} ensure that the shipments between two trucks is occurred at the same period or before docking period of the outbound truck. Finally, as presented in the main model, the capacity constraints are imposed to this step of the solution approach. 

The decomposition of the main problem into two sub-problems decrease significantly the complexity of the problem. But a large number of variables and constraints still remain in the first sub-problem. For this reason, an adaptive heuristic method is proposed which could be employed instead of the mathematical model.   

\subsection{Adaptive heuristic}

A heuristic algorithm is developed to solve the first sub-problem and provide the good feasible solutions. The obtained feasible solution is employed as a starting point for the solver to help that in converging more quickly to the optimal solution. The heuristic algorithm aims, on the one hand, to assign the outbound trucks to the destinations and periods, and on the other hand, to fill up the trucks according to the demands of their assigned destinations. The steps of the heuristic are summarized in \autoref{Fig-Heuristic} and presented as follows:

\begin{figure}

\centering
\includegraphics[width=0.5\textwidth]{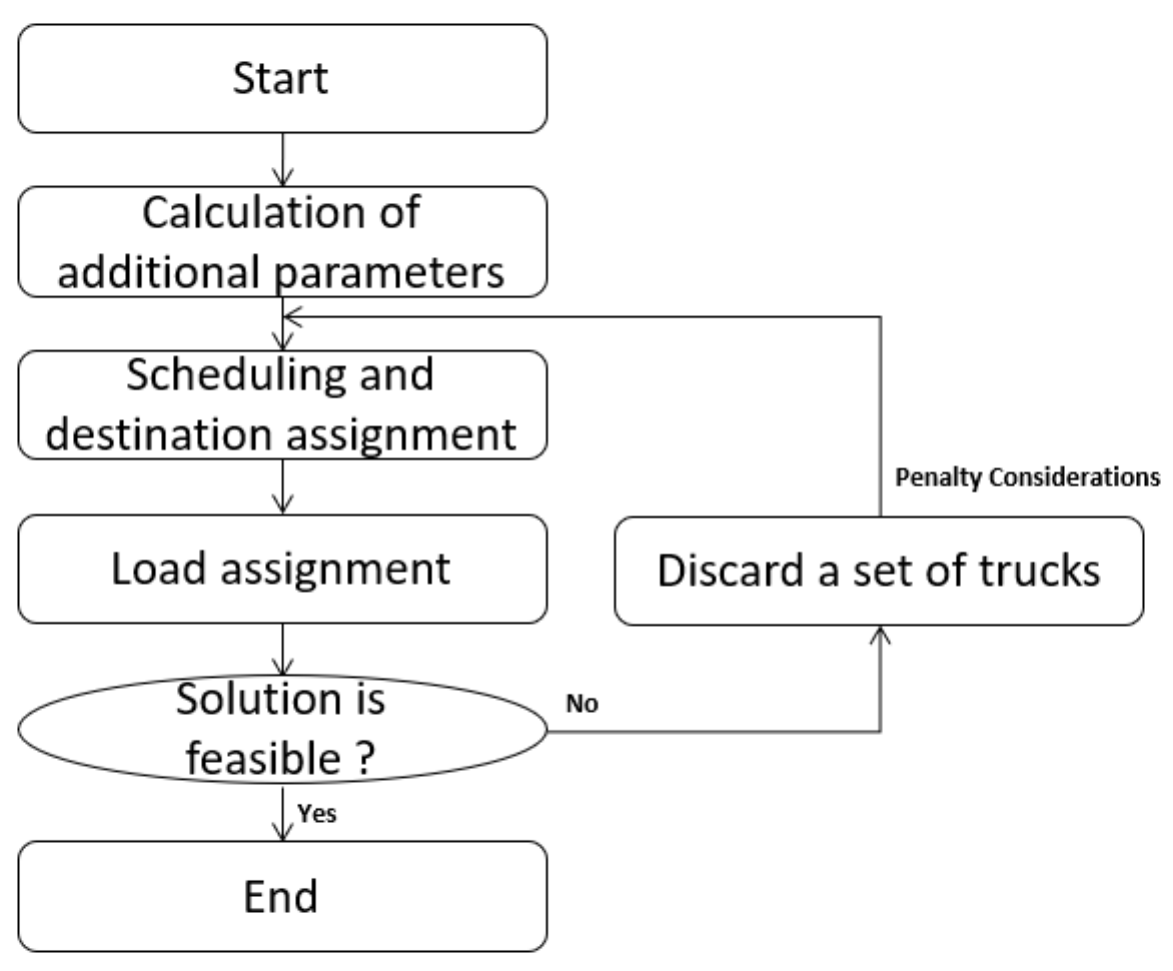}
\caption{Heuristic algorithm}
\label{Fig-Heuristic}
\end{figure}

\subsubsection{Pre-processing and calculations}

The first step of the algorithm concerns with the definition and the calculation of three auxiliary parameters which are obtained from the parameters of the main model. These parameters are needed for the next steps and they are
listed as follows:

\begin{itemize}
\item The total demand of every destination at each time period : $RP_{dt}=\sum_{p\in P}R_{pdt};\;\forall d\in D,\:t\in T;$ 
\item The number of trucks needed for every destination in each time period
: $ZP_{dt}=Ciel(RP_{dt}/C);\;\forall d\in D,\:t\in T;$
\item The total number of trucks needed in each time period : $ZT_{t}=\sum_{p\in P}ZP_{dt};\;\forall t\in T;$
\end{itemize}

In addition, the detailed demand of each destination ($R_{pdt}$) and the total number of outbound trucks $n$ are used to schedule and assign the outbound trucks in the next steps. 

\subsubsection{Scheduling and assignment of outbound trucks to the destinations}

The algorithm starts by assigning trucks to the destinations that need products at the first time period in ascending order of the truck numbers and destination numbers. Considering that the trucks are homogeneous and they have no restriction or priority, they can be assigned to any destination at any time period. Note that, a truck could not be assigned to more than one destination, but a destination could need more that one truck in a time period. In the same manner, the trucks are attributed to the destinations over next periods. In this way, a number of trucks equal to ``$ZT_{t}$'', are programmed to be docked in period $t$. This step of the algorithm is terminated by assigning a needed number of trucks to the last required destination in the last planning period.

Actually, by having the total number of pallets which is needed by each destination at every time period ($RP_{dt}$), the number of trucks needed to deliver these pallets ($ZP_{dt}$) is determined. The first truck (group of trucks) is attributed to the first destination in the first period (if this destination needs a shipment). The next truck (group of trucks) is assigned to the next destination which needs a shipment in the same planning period and so on until the demand of the last destination. The same procedure is employed to assign the remained trucks to the destinations demands at the next periods, until any demand is attributed to a truck. At the end of this step, the destination of each outbound truck ($q_{jd}$), the time period in which the trucks are docked ($y_{jt}$) and the number of pallets placed in each outbound truck ($\sum_{p\in P}x_{jp}$) are determined, but the type of pallets (products) is unknown. The procedure of this step is summarized in algorithm \ref{Alg-S1}.

\begin{algorithm}

\caption{Scheduling and assignment of outbound trucks to the destinations}
\label{Alg-S1}

{\footnotesize{}}%
\begin{tabular}{lc}
{\footnotesize{}$j=0;$} & \tabularnewline
{\footnotesize{}$for(t\in T)\{for(d\in D)\{$ } & \tabularnewline
{\footnotesize{}$\quad BP_{dt}=RP_{dt};$} & {\footnotesize{}//initializing the auxiliary variables}\tabularnewline
{\footnotesize{}$\quad while(BP_{dt}>0)\{$} & \tabularnewline
{\footnotesize{}$\qquad y_{jt}=1;\;d_{jd}=1;\;LJ_{j}=min\{C,BP_{dt}\}$;} & {\footnotesize{}//updating the decision variables }\tabularnewline
{\footnotesize{}$\qquad BP_{dt}=BP_{dt}-LJ_{j};\;j=j+1;\}$ } & {\footnotesize{}//updating the auxiliary variables}\tabularnewline
{\footnotesize{}$\}\}$} & \tabularnewline
\end{tabular}{\footnotesize\par}

\end{algorithm}

\subsubsection{Assignment of loads to the outbound trucks}

In the previous step, according to the total demand of the destinations in each period, the total number of pallets attributed to each truck is determined. By considering the detailed destination demands ($R_{pdt}$), in this step, the type of pallets which must be placed in each truck is decided. The number of pallets of each type in a truck is related to the demands of the destination to which the truck is assigned. In the studied JIT distribution system, each truck is filled with the needed pallets at each period. Hence, a group of trucks may be assigned to the same destination in different planning periods, but the load of each one depends on the period to which it is assigned and the demand of the attributed destination in that period ($R_{pdt}$).

The algorithm starts with the first planning period, first destination and first product type. It evaluates if the related parameter $R$ is positive or not. If $R>0$, then, it finds the first truck which is assigned to this destination in this planning period. If the truck has a free capacity $B_{j}>0$, a load of the considered product is attributed to that $x_{jp}=min\{B_{j}\:\&\:R_{pdt}\}$. Subsequently, the remaining capacity of the truck and destination demand are revised. If the truck capacity is filled before the destination demand, then, the next truck assigned to the destination in current planning period must be found and pooled similarly. Otherwise, the algorithm passes to the next product, ordered by the same destination in the current planning period. This procedure continues until the assigned trucks are filled with all needed products for theirs designated destinations during the first period. For the next periods, the pallets of products (loads or shipments) are assigned to the trucks in the same manner. The details of the load of each truck, i.e. number of pallets of each product ($x_{jp}$) are determined in this step. Algorithm \ref{Alg-S2} shows the general procedure of this step in summary.

\begin{algorithm}
\caption{Assignment of loads to the outbound trucks}
\label{Alg-S2}

{\footnotesize{}}%
\begin{tabular}{lc}
{\footnotesize{}$for(t\in T)\{for(d\in D)\{for(p\in P)\{$} & \tabularnewline
{\footnotesize{}$\quad B_{pdt}=R_{pdt};$} & {\footnotesize{}//initializing the auxiliary variables}\tabularnewline
{\footnotesize{}$\quad for(j\in J)\{$} & \tabularnewline
{\footnotesize{}$\qquad if(y_{jt}>0\:and\:q_{jd}>0\:and\:B_{pdt}>0\:and\:BL_{j}>0)\{$} & \tabularnewline
{\footnotesize{}$\quad\qquad x_{jp}=min\{B_{pdt},BL_{j}\};\;d_{jd}=1;\;LJ_{j}=min\{C,BP_{dt}\}$;} & {\footnotesize{}//updating the decision variables }\tabularnewline
{\footnotesize{}$\quad\qquad BL_{j}=BL_{j}-x_{jp};\;B_{pdt}=B_{pdt}-x_{jp};\}$} & {\footnotesize{}//updating the auxiliary variables}\tabularnewline
{\footnotesize{}$\}\}$\}} & \tabularnewline
\end{tabular}{\footnotesize\par}
\end{algorithm}

\subsubsection{Feasibility evaluation}

Three previous steps prepare a solution that satisfies all destinations demands at their needed periods. This solution may be unfeasible because of the limited number of outbound doors. Actually, the number of trucks that can be docked to the cross-dock simultaneously (in the same period) must not be more than the number of doors. Otherwise, a number of trucks must be ignored in the violating periods. That signifies the non-delivered loads. Hence, the algorithm chooses the trucks with minimum loads to be removed. 

Mathematically, the number of trucks needed in each period ($ZT_{t}$) was calculated in the second step of the algorithm. If $ZT_{t}>OD$, then, among the pooled trucks in the current period, algorithm finds the truck with minimum load ($min_{j\in J^{t}}\{L_{j}\}$); Where $J^{t}$ is a subset of $J$ that signifies the trucks which are docked in period $t$, and $L_{j}=\sum_{p\in P}x_{jp}$is the total pallets placed in truck $j$. This cycle is repeated until the number of assigned truck is equal to the number of door for every period ($ZT_{t}=OD$). Consequently, a feasible solution is obtained that satisfies the maximum number of destinations demands. After removing some trucks of the periods violating the outbound door constraint, the lists of the scheduled trucks in each period ($y_{jt}$) and the assigned truck to each destination ($q_{jd}$) are revised.

\section{Computational results}\label{results}

In order to evaluate the solution approaches, the main (integrated) model is analyzed in the first step. Then, the symmetry breaking constraints are added to the model to investigate the impact of these constraints on the efficiency of the model. In the second step, the matheuristic method are evaluated by solving various instances and that is compared with the best scenario of the first step (integrated model with symmetry breaking constraints). In view of the complexity of the first step of matheuristic to be solved by CPLEX, an adaptive heuristic algorithm is employed to solve that. In this way, two scenarios are proposed. In the first one, the solution of the heuristic method is injected as a primal (initial) solution for the first step of matheuristic, then that is solved by Cplex. In the second scenario, the heuristic algorithm is employed as a solution method to the first step and the solutions are directly transferred to the second step of matheuristic. These two scenarios are analyzed and are compared with each other and with the classical matheuristic in which both steps are solved by Cplex. Note that, in the numerical experiments, a working day is decomposed to $5$ periods of $2$ hours. In a multi-commodity logic, the products which must be distributed across the studied cross-dock are classified into five categories. Hence, two parameters are constant in all of the studied instances ($r=5;\;k=5$). 

\subsection{Impact of the symmetry braking constraints}

To evaluate the impact of the symmetry braking constraints, a sum of $23$ small and medium size instances are solved with and without these constraints. A time limit of $10800$ seconds is imposed to the solver. The results of the experiments are presented in table \ref{tab:SB Impact}.

\begin{table}
\caption{\label{tab:SB Impact}Impact of the Symmetry Breaking constraints}
\scriptsize

$r:Periods=5$

$k:Product\;types=5$

\begin{adjustbox}{max width=\textwidth}

\begin{tabular}{cccccccccccc}
\hline 
\multirow{2}{*}{DB} & \multirow{2}{*}{Ins} & \multirow{2}{*}{IT-OT-D-ID-OD} & \multicolumn{4}{c}{With Symmetry Breaking} &  & \multicolumn{4}{c}{Without Symmetry Breaking}\tabularnewline
\cline{4-7} \cline{5-7} \cline{6-7} \cline{7-7} \cline{9-12} \cline{10-12} \cline{11-12} \cline{12-12} 
 &  &  & CPU & Nodes & BFS & Status &  & CPU & Nodes & BFS & Status\tabularnewline
\hline 
\multirow{3}{*}{1} & 1.1 & 5 - 8 - 3 - 3 - 3 & 1 & 7857 & 1304{*} & Opt &  & 8 & 14236 & 1304{*} & Opt\tabularnewline
 & 1.2 & 5 - 7 - 3 - 3 - 3 & <1 & 1037 & 13{*} & Opt &  & 31 & 72160 & 13{*} & Opt\tabularnewline
 & 1.3 & 5 - 9 - 3 - 3 - 3 & 10 & 6903 & 1105{*} & Opt &  & 35 & 91452 & 1105{*} & Opt\tabularnewline
\hline 
\multirow{4}{*}{2} & 2.1 & 6 - 11 - 4 - 4 - 4 & 13 & 9507 & 14{*} & Opt &  & >10800 & 1.54e+6 & 14 & 40.05\%\tabularnewline
 & 2.2 & 6 - 9 - 4 - 4 - 4 & 7 & 5446 & 1205{*} & Opt &  & 14 & 16304 & 1205{*} & Opt\tabularnewline
 & 2.3 & 6 - 10 - 4 - 4 - 4 & 10 & 6291 & 214{*} & Opt &  & 12 & 16579 & 214{*} & Opt\tabularnewline
 & 2.4 & 6 - 11 - 4 - 4 - 4 & 62 & 121824 & 1510{*} & Opt &  & >10800 & 336159 & 1510 & 0.66\%\tabularnewline
\hline 
\multirow{4}{*}{3} & 3.1 & 7 - 11 - 4 - 4 - 4 & 19 & 5930 & 1905{*} & Opt &  & 415 & 346174 & 1905{*} & Opt\tabularnewline
 & 3.2 & 7 - 12 - 4 - 4 - 4 & 44 & 42589 & 1608{*} & Opt &  & >10800 & 397552 & 1608 & 0.37\%\tabularnewline
 & 3.3 & 7 - 10 - 4 - 4 - 4 & 14 & 7192 & 614{*} & Opt &  & 129 & 221713 & 614{*} & Opt\tabularnewline
 & 3.4 & 7 - 9 - 4 - 4 - 4 & 5 & 6117 & 1103{*} & Opt &  & 6 & 6919 & 1103{*} & Opt\tabularnewline
\hline 
\multirow{4}{*}{4} & 4.1 & 8 - 12 - 5 - 5 - 5 & 638 & 1.89e+6 & 3208{*} & Opt &  & >10800 & 960857 & 3208 & 0.25\%\tabularnewline
 & 4.2 & 8 - 14 - 5 - 5 - 5 & 59 & 41205 & 410{*} & Opt &  & >10800 & 279852 & 410 & 2.03\%\tabularnewline
 & 4.3 & 8 - 12 - 5 - 5 - 5 & 9 & 9230 & 308{*} & Opt &  & >10800 & 594711 & 308 & 2.60\%\tabularnewline
 & 4.4 & 8 - 13 - 5 - 5 - 5 & 13 & 5112 & 218{*} & Opt &  & >10800 & 1.00e+6 & 218 & 7.19\%\tabularnewline
\hline 
\multirow{4}{*}{5} & 5.1 & 9 - 13 - 6 - 6 - 6 & 2152 & 3.58e+6 & 1604{*} & Opt &  & >10800 & 308970 & 1604 & 19.89\%\tabularnewline
 & 5.2 & 9 - 17 - 6 - 6 - 6 & 44 & 42926 & 10{*} & Opt &  & >10800 & 239477 & 10 & 100\%\tabularnewline
 & 5.3 & 9 - 15 - 6 - 6 - 6 & 280 & 50119 & 616{*} & Opt &  & >10800 & 201172 & 616 & 99.68\%\tabularnewline
 & 5.4 & 9- 17 - 6 - 6 - 6 & 139 & 43869 & 514{*} & Opt &  & >10800 & 197311 & 514 & 2.72\%\tabularnewline
\hline 
\multirow{4}{*}{6} & 6.1 & 10 - 18 - 7 - 7 - 7 & 159 & 66644 & 507{*} & Opt &  & >10800 & 193237 & 507 & 100\%\tabularnewline
 & 6.2 & 10 - 20 - 7 - 7 - 7 & >10800 & 1.01e+6 & 1105 & 0.27\% &  & >10800 & 114208 & 1106 & 100\%\tabularnewline
 & 6.3 & 10 - 16 - 7 - 7 - 7 & 3249 & 643603 & 806{*} & Opt &  & >10800 & 168406 & 806 & 62.78\%\tabularnewline
 & 6.4 & 10 - 15 - 7 - 7 - 7 & 1149 & 953575 & 617{*} & Opt &  & >10800 & 283059 & 617 & 2.76\%\tabularnewline
\hline 
\end{tabular}
\end{adjustbox}
\end{table}

In this table, ``DB'' means the data-base, ``Ins'' represents the instances where ``IT'', ``OT'', ``D'', ``ID'' and ``OD'' signify the number of Inbound Trucks, the number of Outbound Trucks, the number of Destinations, the number of Inbound Doors and the number of Outbound Doors respectively. In each data-base, the number of inbound trucks, destinations, inbound doors and outbound doors are fixed, the orders and their needed JIT periods are simulated and the number of outbound trucks are calculated based on the received orders. For the small and medium size instances, IT varies between 5 and 10 trucks, OT is from $7$ to $20$ trucks, and D, ID and OD vary between $3$ and $7$. 

Two solution approaches are compared by their CPU times (in second), best feasible solutions (BFS) and the status of the solutions. If the solution approach attains the optimal solution in the considered time limit, the status of the solution is optimal (opt), if not, the status column presents the gap between lower and upper bounds of the solutions found in the time limit in percentage. 

Note that the penalty of the undelivered pallets is considered to be $100$ for any pallet. As shown in table \ref{tab:SB Impact}, symmetry breaking constraints increase the efficiency of the model significantly. Without SBC, the solver attains the optimal solution for only 7 instances, whereas with SBC the solver reach to the optimal solution for 22 instances among 23, within the defined time limit. Even in the cases where both scenarios give the optimal solution, the number of nodes and computational time are decreased with the SBCs. As an example, for instance 3.1, the solver analyzed $340\times10^{3}$ nodes during $415$ second to attain the optimal solution, whereas with SBCs, $5900$ nodes are analyzed in $19$ seconds. Therefore, the model with SBCs overcomes the model without the SBCs over all studied instances. 

Given the exponentially increasing of the computational time by the size of instances, even with SBCs, the generic solver is not suitable for the large size instances. For this reason, a matheuristic method is proposed and evaluated in the next section. 

\subsection{Numerical results of matheuristic method}

A set of $39$ instances of different sizes (small, medium and large instances) are employed to evaluate the matheuristic method. Three essential parameters of the problem, which determine the size of the instances, are as follows: the number of inbound trucks, the number of outbound trucks and the number of destinations. The number of inbound trucks varies from $5$ to $80$, the variation range for the outbound trucks is between $7$ and $110$, where the number of destinations varies from $3$ to $40$ to generate the instances. In order to evaluate the matheuristic, the result of the best scenario of the previous analysis (Cplex with SBC) is selected to be compared with the matheuristic results. The proposed two-step matheuristic is coded by c++ language and both steps are solved by ILOG CPLEX. The time limit considered for each step is $5400$ seconds and for the complete solution procedure is $10,800$ seconds. The numerical results of the matheuristic are presented in table \ref{tab:MathResults} and compared with the results of CPLEX.

\begin{table}
\caption{\label{tab:MathResults}Results of the matheuristic method}

\scriptsize
\begin{adjustbox}{max width=\textwidth}

\begin{tabular}{cccccccccccccccc}
\hline 
\multirow{3}{*}{Ins} & \multirow{3}{*}{IT-OT-D} & \multicolumn{4}{c}{Cplex with Symmetry Breaking} &  & \multicolumn{9}{c}{2-step Matheuristic method}\tabularnewline
\cline{8-16} \cline{9-16} \cline{10-16} \cline{11-16} \cline{12-16} \cline{13-16} \cline{14-16} \cline{15-16} \cline{16-16} 
 &  & \multicolumn{4}{c}{} &  & \multicolumn{3}{c}{Step1 (Cplex)} &  & \multicolumn{3}{c}{Step2 (Cplex)} & \multirow{2}{*}{BFS} & \multirow{2}{*}{E$_{max}$}\tabularnewline
\cline{3-6} \cline{4-6} \cline{5-6} \cline{6-6} \cline{8-10} \cline{9-10} \cline{10-10} \cline{12-14} \cline{13-14} \cline{14-14} 
 &  & CPU & Nodes & BFS & Status &  & CPU & Nodes & Status &  & CPU & Node & Status &  & \tabularnewline
\hline 
1.1 & 5 - 8 - 3 & 1 & 7857 & 1304{*} & Opt &  & <1 & 33 & Opt &  & <1 & 0 & Opt & 1309 & 0.38\%{*}\tabularnewline
1.2 & 5 - 7 - 3 & <1 & 1037 & 13{*} & Opt &  & <1 & 0 & Opt &  & <1 & 0 & Opt & 13{*} & Opt\tabularnewline
1.3 & 5 - 9 - 3 & 10 & 6903 & 1105{*} & Opt &  & <1 & 3 & Opt &  & <1 & 0 & Opt & 1106 & 0.09\%{*}\tabularnewline
2.1 & 6 - 11 - 4 & 13 & 9507 & 14{*} & Opt &  & <1 & 0 & Opt &  & <1 & 544 & Opt & 14{*} & Opt\tabularnewline
2.2 & 6 - 9 - 4 & 7 & 5446 & 1205{*} & Opt &  & <1 & 5 & Opt &  & <1 & 0 & Opt & 1206 & 0.08\%{*}\tabularnewline
2.3 & 6 - 10 - 4 & 10 & 6291 & 214{*} & Opt &  & <1 & 0 & Opt &  & <1 & 35 & Opt & 214{*} & Opt\tabularnewline
2.4 & 6 - 11 - 4 & 62 & 121824 & 1510{*} & Opt &  & <1 & 270 & Opt &  & <1 & 529 & Opt & 1512 & 0.13\%{*}\tabularnewline
3.1 & 7 - 11 - 4 & 19 & 5930 & 1905{*} & Opt &  & <1 & 45 & Opt &  & <1 & 183 & Opt & 1908 & 0.15\%{*}\tabularnewline
3.2 & 7 - 12 - 4 & 44 & 42589 & 1608{*} & Opt &  & <1 & 0 & Opt &  & <1 & 4114 & Opt & 1609 & 0.06\%{*}\tabularnewline
3.3 & 7 - 10 - 4 & 14 & 7192 & 614{*} & Opt &  & <1 & 0 & Opt &  & <1 & 369 & Opt & 618 & 0.64\%{*}\tabularnewline
3.4 & 7 - 9 - 4 & 5 & 6117 & 1103{*} & Opt &  & <1 & 33 & Opt &  & <1 & 164 & Opt & 1105 & 0.18\%{*}\tabularnewline
4.1 & 8 - 12 - 5 & 638 & 1.89e+6 & 3208{*} & Opt &  & <1 & 0 & Opt &  & 1.22 & 4661 & Opt & 3209 & 0.03\%{*}\tabularnewline
4.2 & 8 - 14 - 5 & 59 & 41205 & 410{*} & Opt &  & <1 & 0 & Opt &  & 1.13 & 3176 & Opt & 410{*} & Opt\tabularnewline
4.3 & 8 - 12 - 5 & 9 & 9230 & 308{*} & Opt &  & <1 & 0 & Opt &  & <1 & 208 & Opt & 308{*} & Opt\tabularnewline
4.4 & 8 - 13 - 5 & 13 & 5112 & 218{*} & Opt &  & <1 & 461 & Opt &  & <1 & 81 & Opt & 222 & 1.83\%{*}\tabularnewline
5.1 & 9 - 13 - 6 & 2152 & 3.58e+6 & 1604{*} & Opt &  & <1 & 153 & Opt &  & 2.33 & 4109 & Opt & 1604 & 0.16\%{*}\tabularnewline
5.2 & 9 - 17 - 6 & 44 & 42926 & 10{*} & Opt &  & <1 & 0 & Opt &  & 2.23 & 6025 & Opt & 12 & 20\%{*}\tabularnewline
5.3 & 9 - 15 - 6 & 280 & 50119 & 616{*} & Opt &  & <1 & 0 & Opt &  & 3.44 & 4106 & Opt & 619 & 0.48\%{*}\tabularnewline
5.4 & 9 - 17 - 6 & 139 & 43869 & 514{*} & Opt &  & <1 & 0 & Opt &  & 1.03 & 358 & Opt & 516 & 0.38\%{*}\tabularnewline
6.1 & 10 - 18 - 7 & 159 & 66644 & 507{*} & Opt &  & <1 & 0 & Opt &  & 6.92 & 7274 & Opt & 508 & 0.19\%{*}\tabularnewline
6.2 & 10 - 20 - 7 & >10800 & 1.01e+6 & 1105 & 0.27\% &  & 2.22 & 560 & Opt &  & 4.81 & 5663 & Opt & 1106 & 0.09\%\tabularnewline
6.3 & 10 - 16 - 7 & 3249 & 643603 & 806{*} & Opt &  & <1 & 0 & Opt &  & 2.89 & 4238 & Opt & 806{*} & Opt\tabularnewline
6.4 & 10 - 15 - 7 & 1149 & 953575 & 617{*} & Opt &  & <1 & 0 & Opt &  & <1 & 193 & Opt & 620 & 0.48\%{*}\tabularnewline
7.1 & 12 - 20 - 8 & >10800 & 3.49e+6 & 1206 & 0.43\% &  & 3.42 & 0 & Opt &  & 18.33 & 14046 & Opt & 1206{*} & 0\%\tabularnewline
8.1 & 15 - 22 - 9 & >10800 & 710370 & 930 & 1.89\% &  & 1.47 & 196 & Opt &  & 4.20 & 3966 & Opt & 931 & 0.10\%\tabularnewline
9.1 & 17 - 24 - 10 & >10800 & 136119 & 19 & 64.18\% &  & 1.30 & 0 & Opt &  & 449.81 & 110571 & Opt & 20 & 5.2\%\tabularnewline
10.1 & 20 - 24 - 10 & >10800 & 309978 & 1315 & 1.14\% &  & 2.55 & 1046 & Opt &  & >10800 & 303973 & 0.48\% & 1316 & 0.07\%\tabularnewline
11.1 & 22 - 37 - 12 & >10800 & 38353 & 807 & 0.87\% &  & 26.06 & 0 & Opt &  & >10800 & 451043 & 1.07\% & 809 & 0.24\%\tabularnewline
12.1 & 25 - 38 - 15 & >10800 & 12928 & 2818 & 100\% &  & 80.77 & 5522 & Opt &  & >10800 & 139045 & 0.20\% & \textbf{2714} & \textbf{-3.69\%}\tabularnewline
13.1 & 27 - 45 - 15 & >10800 & 14853 & 2919 & 100\% &  & 198.31 & 4384 & Opt &  & >10800 & 1.18e+6 & 0.17\% & \textbf{1803} & \textbf{-38.23\%}\tabularnewline
14.1 & 30 - 49 - 15 & >10800 & 14254 & 12136 & 95.06\% &  & >5400 & 121143 & 34.04\% &  & >5400 & 229241 & 0.02\% & \textbf{5701} & \textbf{-53.02\%}\tabularnewline
15.1 & 32 - 54 - 18 & >10800 & 15910 & 11601 & 100\% &  & 1039.34 & 19571 & Opt &  & >10800 & 284247 & 0.17\% & \textbf{1803} & \textbf{-84.45\%}\tabularnewline
16.1 & 35 - 56 - 19 & >10800 & 7772 & 10520 & 100\% &  & 2823.49 & 6498 & Opt &  & >10800 & 857371 & 100\% & \textbf{6} & \textbf{-99.94\%}\tabularnewline
17.1 & 40 - 63 - 20 & >10800 & - & - & - &  & 3497.91 & 33876 & Opt &  & >10800 & 427966 & 0.16\% & \textbf{3105} & -\tabularnewline
18.1 & 45 - 67 - 23 & >10800 & - & - & - &  & >5400 & 15661 & 92.16\% &  & >5400 & 321555 & 0.14\% & \textbf{5107} & -\tabularnewline
19.1 & 50 - 79 - 25 & >10800 & - & - & - &  & >5400 & 5764 & 100\% &  & >5400 & 334438 & 0.08\% & \textbf{5104} & -\tabularnewline
20.1 & 60 - 94 - 30 & >10800 & - & - & - &  & >5400 & 3955 & 100\% &  & 1451.55 & 188938 & Opt & \textbf{10701} & -\tabularnewline
21.1 & 70 - 109 - 35 & >10800 & - & - & - &  & >5400 & 112 & 100\% &  & 1.34 & 0 & Opt & \textbf{118910} & -\tabularnewline
22.1 & 80 - 115 - 40 & >10800 & - & - & - &  & >5400 & 1385 & 100\% &  & 20.88 & 0 & Opt & \textbf{74303} & -\tabularnewline
\hline 
\end{tabular}
\end{adjustbox}
\end{table}

As shown in the table, by employing the SBCs and Cplex solver, the optimal solutions are found for 22 small and medium instances in the limited time available (shown by {*}). For 5 medium instances, the solver finds the efficient lower and upper bounds for the solutions (Gaps are less than 2\%). But for 12 large instances, Cplex solver whether is not able to find two bounds (6 instances) or is not able to finds an acceptable lower bound (Gaps are greater than 64\%). The BFS column of the matheuristic method shows that among 22 instances for which the optimal solutions are known, 6 instances are solved
to optimality (shown by {*}). The $E_{max}$ column of the table shows the status of the solutions obtained by matheuristic. If the solution is not optimal, that shows the gaps between the found solution and the optimal solution. Note that, if the optimal solution is available, the gaps are shown together with a {*}, if not, the $E_{max}$ presents the difference between the solutions found by CPLEX and matheuristic method. For 11 large instances, the matheuristic method attains a solution better than Cplex (shown by bold), where for the other instances the gaps are mostly less than 1\% expect 3 instances (Instance 4.4 : 1.83\%, Instance 9.1 : 5.2\% and Instance 5.2 : 20\%). 

Table \ref{tab:MathResults} shows that by using the matheuristic method the complexity of the problem decreases greatly but the quality of the solutions do not decrease significantly. The number of analyzed nodes and the CPU times illustrate the efficiency of the proposed solution method. As an example, for instance 6.3, Cplex analyzed more than 600 thousands nodes during more than 3 thousands seconds to find the optimal solution whereas matheuristic analyze about 4 thousands nodes within 3 seconds to reach the same solution as CPLEX (optimal solution).

As the first step of the proposed matheuristic has a significant impact on the quality of the solution, this method is not very suitable where the solution of the first step is not reliable. For the five largest instances, the solver can not attains an optimal solution or a solution with a small gap in the first step. Hence, the proposed classical matheuristic is not a qualified approach to solve these sizes of the problem. For this reason, a heuristic algorithm is presented to solve the first step of the matheuristic. 

\subsection{Warmed-up matheuristic vs hybrid matheuristic}

The heuristic algorithm, which is proposed for the first step of matheuristic method, provides two scenarios to solve the problem. The first scenario is to apply the heuristic method as an alternative to solve the first step instead of CPLEX. In the second scenario, the solutions obtained by the heuristic method are used as the primal solution to be injected to CPLEX. In this way, the first step is solved by CPLEX, but from a good starting solution obtained by heuristic. These two scenarios are compared to the classical matheuristic method (CPLEX-CPLEX) by 39 studied instances. The results of the classical matheuristic method is presented in table \ref{tab:MathResults}. Table \ref{tab:WU-Hybride} shows the computational results of two aforementioned alternatives for the first step (1. CPLEX with a primal solution - CPLEX, 2. Heuristic - CPLEX).

\begin{flushleft}
\begin{table}
\caption{\label{tab:WU-Hybride}Results of matheuristic with Warm-Up and hybrid
matheuristic}

\scriptsize
\begin{adjustbox}{max width=\textwidth}
\begin{raggedright}
\begin{tabular}{cccccccccccccccccccc}
\hline 
\multirow{3}{*}{Ins} & \multicolumn{9}{c}{Matheuristic with Warm-Up} &  & \multicolumn{9}{c}{Hybrid matheuristic}\tabularnewline
\cline{2-10} \cline{3-10} \cline{4-10} \cline{5-10} \cline{6-10} \cline{7-10} \cline{8-10} \cline{9-10} \cline{10-10} \cline{12-20} \cline{13-20} \cline{14-20} \cline{15-20} \cline{16-20} \cline{17-20} \cline{18-20} \cline{19-20} \cline{20-20} 
 & \multicolumn{3}{c}{Step1 (Warmed Up Cplex)} &  & \multicolumn{3}{c}{Step2 (Cplex)} & \multirow{2}{*}{BFS} & \multirow{2}{*}{E$_{max}$} &  & \multicolumn{3}{c}{Step1 (Heuristic)} &  & \multicolumn{3}{c}{Step2 (Cplex)} & \multirow{2}{*}{BFS} & \multirow{2}{*}{E$_{max}$}\tabularnewline
\cline{2-4} \cline{3-4} \cline{4-4} \cline{6-8} \cline{7-8} \cline{8-8} \cline{12-14} \cline{13-14} \cline{14-14} \cline{16-18} \cline{17-18} \cline{18-18} 
 & CPU & Nodes & Status &  & CPU & Node & Status &  &  &  & CPU & Sol & E$_{m}$ &  & CPU & Node & Status &  & \tabularnewline
\hline 
1.1 & <1 & 0 & Opt &  & <1 & 0 & Opt & 1304{*} & Opt &  & <1 & 13{*} & Opt &  & <1 & 0 & Opt & 1304{*} & Opt\tabularnewline
1.2 & <1 & 0 & Opt &  & <1 & 0 & Opt & 13{*} & Opt &  & <1 & 0{*} & Opt &  & <1 & 0 & Opt & 13{*} & Opt\tabularnewline
1.3 & <1 & 3 & Opt &  & <1 & 0 & Opt & 1105{*} & Opt &  & <1 & 11{*} & Opt &  & <1 & 0 & Opt & 1105{*} & Opt\tabularnewline
2.1 & <1 & 0 & Opt &  & <1 & 536 & Opt & 14{*} & Opt &  & <1 & 0{*} & Opt &  & <1 & 536 & Opt & 14{*} & Opt\tabularnewline
2.2 & <1 & 35 & Opt &  & <1 & 0 & Opt & 1206 & 0.08\%{*} &  & <1 & 12{*} & Opt &  & <1 & 0 & Opt & 1206 & 0.08\%{*}\tabularnewline
2.3 & <1 & 0 & Opt &  & <1 & 48 & Opt & 214{*} & Opt &  & <1 & 2{*} & Opt &  & <1 & 48 & Opt & 214{*} & Opt\tabularnewline
2.4 & <1 & 202 & Opt &  & <1 & 594 & Opt & 1511 & 0.06\%{*} &  & <1 & 15{*} & Opt &  & <1 & 594 & Opt & 1511 & 0.06\%{*}\tabularnewline
3.1 & <1 & 311 & Opt &  & <1 & 240 & Opt & 1906 & 0.05\%{*} &  & <1 & 19{*} & Opt &  & <1 & 240 & Opt & 1906 & 0.05\%{*}\tabularnewline
3.2 & <1 & 788 & Opt &  & <1 & 1576 & Opt & 1609 & 0.06\%{*} &  & <1 & 16{*} & Opt &  & <1 & 1576 & Opt & 1609 & 0.06\%{*}\tabularnewline
3.3 & <1 & 0 & Opt &  & <1 & 387 & Opt & 617 & 0.48\%{*} &  & <1 & 6{*} & Opt &  & <1 & 387 & Opt & 617 & 0.48\%{*}\tabularnewline
3.4 & <1 & 40 & Opt &  & <1 & 66 & Opt & 1104 & 0.09\%{*} &  & <1 & 11{*} & Opt &  & <1 & 66 & Opt & 1104 & 0.09\%{*}\tabularnewline
4.1 & <1 & 5 & Opt &  & <1 & 1510 & Opt & 3209 & 0.03\%{*} &  & <1 & 32{*} & Opt &  & <1 & 1510 & Opt & 3209 & 0.03\%{*}\tabularnewline
4.2 & <1 & 23 & Opt &  & 1.64 & 3078 & Opt & 410{*} & Opt &  & <1 & 4{*} & Opt &  & 1.48 & 3078 & Opt & 410{*} & Opt\tabularnewline
4.3 & <1 & 0 & Opt &  & <1 & 391 & Opt & 308{*} & Opt &  & <1 & 3{*} & Opt &  & <1 & 391 & Opt & 308{*} & Opt\tabularnewline
4.4 & <1 & 1187 & Opt &  & <1 & 788 & Opt & 220 & 0.91\%{*} &  & <1 & 2{*} & Opt &  & <1 & 788 & Opt & 220 & 0.91\%{*}\tabularnewline
5.1 & 1.00 & 3441 & Opt &  & 1.70 & 1696 & Opt & 1605 & 0.06\%{*} &  & <1 & 16{*} & Opt &  & 1.48 & 1696 & Opt & 1605 & 0.06\%{*}\tabularnewline
5.2 & <1 & 0 & Opt &  & 1.77 & 3667 & Opt & 11 & 10\%{*} &  & <1 & 0{*} & Opt &  & 1.67 & 3667 & Opt & 11 & 10\%{*}\tabularnewline
5.3 & <1 & 184 & Opt &  & 3.00 & 1534 & Opt & 618 & 0.32\%{*} &  & <1 & 6{*} & Opt &  & 2.78 & 1534 & Opt & 618 & 0.32\%{*}\tabularnewline
5.4 & 2.89 & 4070 & Opt &  & <1 & 245 & Opt & 515 & 0.19\%{*} &  & <1 & 5{*} & Opt &  & <1 & 245 & Opt & 515 & 0.19\%{*}\tabularnewline
6.1 & 2.81 & 5585 & Opt &  & 9.16 & 7054 & Opt & 509 & 0.39\%{*} &  & <1 & 5{*} & Opt &  & 9.27 & 7054 & Opt & 509 & 0.39\%{*}\tabularnewline
6.2 & 2.19 & 1876 & Opt &  & 2.67 & 1539 & Opt & 1105 & 0\% &  & <1 & 11{*} & Opt &  & 2.77 & 1539 & Opt & 1105 & 0\%\tabularnewline
6.3 & <1 & 758 & Opt &  & 4.00 & 3014 & Opt & 806{*} & Opt &  & <1 & 8{*} & Opt &  & 3.91 & 3014 & Opt & 806{*} & Opt\tabularnewline
6.4 & <1 & 0 & Opt &  & <1 & 934 & Opt & 618 & 0.16\%{*} &  & <1 & 6{*} & Opt &  & <1 & 934 & Opt & 618 & 0.16\%{*}\tabularnewline
7.1 & 29.86 & 21443 & Opt &  & 14.83 & 22344 & Opt & 1206 & 0\% &  & <1 & 12{*} & Opt &  & 24.44 & 22344 & Opt & 1206 & 0\%\tabularnewline
8.1 & <1 & 2135 & Opt &  & 2.41 & 1440 & Opt & 931 & 0.10\% &  & <1 & 9{*} & Opt &  & 2.53 & 1440 & Opt & 931 & 0.10\%\tabularnewline
9.1 & <1 & 0 & Opt &  & 258.86 & 87012 & Opt & 20 & 5.26\% &  & <1 & 0{*} & Opt &  & 257.09 & 87012 & Opt & 20 & 5.26\%\tabularnewline
10.1 & >5400 & 415941 & 100\% &  & >5400 & 196686 & 0.57\% & 1315 & 0\% &  & <1 & 13{*} & Opt &  & >10800 & 415798 & 0.49\% & 1315 & 0\%\tabularnewline
11.1 & >5400 & 536493 & 72.84\% &  & >5400 & 163957 & 0.99\% & 808 & 0.12\% &  & <1 & 8{*} & Opt &  & >10800 & 314531 & 0.82\% & 807 & 0\%\tabularnewline
12.1 & 71.61 & 10483 & Opt &  & >10800 & 254186 & 0.13\% & 2713 & -3.72\% &  & <1 & 27{*} & Opt &  & >10800 & 246762 & 0.13\% & 2713 & -3.72\%\tabularnewline
13.1 & >5400 & 320026 & 16.67\% &  & >5400 & 347553 & 0.22\% & 1804 & -38.19\% &  & <1 & 18{*} & Opt &  & >10800 & 377520 & 0.22\% & 1804 & -38.19\%\tabularnewline
14.1 & >5400 & 171603 & 89.77\% &  & >5400 & 200136 & 0.09\% & 5705 & -52.99\% &  & <1 & 57 & - &  & >10800 & 240900 & 0.09\% & 5705 & -52.99\%\tabularnewline
15.1 & >5400 & 266308 & 100\% &  & >5400 & 157615 & 0.39\% & 1807 & -84.42\% &  & <1 & 18{*} & Opt &  & >10800 & 178702 & 0.39\% & 1807 & -84.42\%\tabularnewline
16.1 & <1 & 0 & Opt &  & 6711.59 & 287102 & Opt & 5 & -99.95\% &  & <1 & 0{*} & Opt &  & 6733.58 & 287039 & Opt & 5 & -99.95\%\tabularnewline
17.1 & >5400 & 36725 & 100\% &  & >5400 & 144137 & 0.26\% & 3108 & - &  & <1 & 31{*} & Opt &  & >10800 & 175891 & 0.26\% & 3108 & -\tabularnewline
18.1 & 2701.16 & 170295 & Opt &  & >10800 & 130821 & 0.22\% & \textbf{5011} & - &  & <1 & 50 & - &  & >10800 & 140860 & 0.22\% & \textbf{5011} & -\tabularnewline
19.1 & >5400 & 34217 & 100\% &  & >5400 & 168164 & 0.12\% & \textbf{4305} & - &  & <1 & 43 & - &  & >10800 & 182593 & 0.12\% & \textbf{4305} & -\tabularnewline
20.1 & >5400 & 10494 & 100\% &  & >5400 & 244870 & 0.12\% & \textbf{2403} & - &  & <1 & 24 & - &  & >10800 & 247599 & 0.12\% & \textbf{2403} & -\tabularnewline
21.1 & >5400 & 1742 & 100\% &  & >5400 & 61024 & 0.06\% & \textbf{7104} & - &  & <1 & 71 & - &  & >10800 & 95676 & 0.04 & \textbf{7103} & -\tabularnewline
22.1 & >5400 & 1652 & 100\% &  & >5400 & 94137 & 0.09\% & \textbf{6506} & - &  & <1 & 65 & - &  & >10800 & 131366 & 0.09\% & \textbf{6506} & -\tabularnewline
\hline 
\end{tabular}
\par\end{raggedright}
\end{adjustbox}
\end{table}
\par\end{flushleft}

By regarding to the best feasible solutions found by these two scenarios and compare them with the classical matheuristic ones presented in table \ref{tab:MathResults}, it can be concluded that both hybrid methods are more efficient than the classical method. The optimal solutions have been found over 8 instances among 22 known optimal solutions. The efficiency of the methods over other instances is evaluated by two parameters : 1. Best found solution (BFS), 2. Deviation from the optimal ($E_{max}$). As shown in table \ref{tab:WU-Hybride}, the BFSs and the deviations of these new solution approaches are smaller than the classical matheuristic ones over most of the instances. But the efficiency of these methods is illustrated by analyzing the five largest instances (18.1 - 22.1) in which the quality of solutions are much more better than matheuristic. 

Furthermore, table \ref{tab:WU-Hybride} shows that the heuristic method, proposed for step 1, attains the optimal solution for all of instances for which the optimal solution is known. The best feasible solutions found by both scenarios and their related deviations from the optimal are the same expect 2 instances (11.1 and 21.1) in which the quality of the solutions found by ``hybrid matheuristic'' are better than the ``warmed up matheuristic'' ones. Hence, the hybrid method overcomes the warmed up matheuristic over all of the studied instances.

As a conclusion, by an efficient mathematical modeling containing SBCs and CPLEX, only the small size instances are solvable (Ins 1.1 to 6.4). For the medium size instances (Ins 7.1 to 17.1), in some instances the hybrid method reaches a better solution and gap and in some other instances the classical method attains the best feasible solutions. Whereas for the large instances (Ins 18.1 to 22.1), the solutions obtained by hybrid approach are much more better than the solutions attained by the classical method. Therefore, the proposed hybrid matheuristic solution approach is very suitable to solve the medium and large sizes instances of the presented truck scheduling
and load assignment problem.  

\section{Conclusions and future research opportunities}\label{conclusion}

A novel operational level optimization problem in the cross-dock planning domain is presented in this study. The studied problem concerns both inbound and outbound scheduling combined with two interrelated (synchronized) assignment problems as : The assignment of products to outbound trucks and the assignment of outbound trucks to destinations. The consideration of the JIT customer demands is another contribution of this research. An integrated mathematical model is presented and strengthened by adding the symmetry breaking constraints. A decomposition based matheuristic and a hybrid matheuristic algorithms are proposed to solve this combinatorial optimization problem. The instances of different sizes are solved by CPLEX. The computational results show that the symmetry breaking constraints can significantly improve the mathematical model. But even by the SBCs, CPLEX is not able to solve the large size instances. Therefore,  a two-step matheuristic based on the decomposition of the formulation is proposed to decompose the integrated model into two smaller mathematical models. The numerical results show that the matheuristic is very efficient for the medium size instances, while the first step of matheuristic remains intractable to be solved by CPLEX. Therefore, an adaptive heuristic is proposed for the first step to form a hybrid matheuristic. The numerical results confirm a significant improvement on the quality of  solutions for the large instances when applying the hybrid matheuristic as the solution approach.

The uncertainty considerations on the arrival times of the trucks to the cross-dock, integrating VRP and split deliveries, adapting an efficient exact algorithm and developing an on line optimization approach could be considered as the future researches.

\bibliographystyle{elsarticle-harv}

\bibliography{Bib}

\end{document}